\documentclass[11pt,a4paper]{article}

\usepackage[labelfont=bf, margin=1em]{caption}

\usepackage[utf8]{inputenc}
\usepackage[T1]{fontenc}
\usepackage[UKenglish]{babel}
\usepackage{lmodern} 

\usepackage{amsmath, amssymb, amsfonts, amsthm}
\usepackage{mathtools} 
\usepackage{physics}   
\usepackage{mathrsfs}  
\usepackage{gensymb}   
\usepackage{textcomp}  
\usepackage[version=4]{mhchem} 
\usepackage{siunitx}   

\usepackage{graphicx}
\usepackage{float}
\usepackage{subcaption}
\usepackage{fullpage}  
\usepackage{appendix}  
\usepackage{setspace}  
\usepackage{booktabs}  
\usepackage{array}
\usepackage{tabularx}  
\usepackage[table]{xcolor} 

\usepackage{enumitem}  
\usepackage{comment}   
\usepackage{listings}  

\usepackage[numbers,sort&compress]{natbib}
\usepackage{url}
\usepackage{hyperref}

\usepackage{marvosym}  
\usepackage{authblk}   

\AtBeginDocument{}
\numberwithin{equation}{section}
\theoremstyle{definition}

\renewcommand{\d}{\mathrm{d}}

\title{\huge{\textbf{Early warning of critical transitions: distinguishing tipping points from Turing destabilizations}}}
\author[1]{Paul A. Sanders}
\author[1]{Robbin Bastiaansen}
\affil[1]{\textit{Mathematical Institute (MI), Institute for Marine and Atmospheric research Utrecht (IMAU) \& Centre for Complex Systems Studies (CCSS),
Utrecht University, 3584 CD Utrecht, The Netherlands}}
\date{\today}

\begin{document}
\maketitle
\begin{abstract}
    Current early warning signs for tipping points often fail to distinguish between catastrophic shifts and less dramatic state changes, such as spatial pattern formation. This paper introduces a novel method that addresses this limitation by providing more information about the type of bifurcation being approached starting from a spatially homogeneous system state. This method relies on estimates of the dispersion relation from noisy spatio-temporal data, which reveals whether the system is approaching a spatially homogeneous (tipping) or spatially heterogeneous (Turing patterning) bifurcation. Using a modified Klausmeier model, we validate this method on synthetic data, exploring its performance under varying conditions including noise properties and distance to bifurcation. We also determine the data requirements for optimal performance. Our results indicate the promise of a new spatial early warning system built on this method to improve predictions of future transitions in many climate subsystems and ecosystems, which is critical for effective conservation and management in a rapidly changing world.
\end{abstract}
\newpage

\section{Introduction}

Nature is changing by the day; climate change and land-use change are putting many ecosystems and climate subsystems under stress. As a result, some systems may be pushed beyond their tipping point -- a critical threshold where even small disturbances can trigger abrupt shifts to alternative states with different functioning. Such tipping points have been identified in many systems~\cite{armstrong2022exceeding, lenton2008tipping, scheffer2001catastrophic, holling1973resilience, drijfhout2015catalogue, terpstra2025assessment}. For example, in the Amazon rainforest where they might induce transitions to savanna-like ecosystems~\cite{Flores2024,Drueke2023}, or in ice sheets where crossing of tipping points indicates they disappear, which might lead to a significant rise in sea level and disrupt coastal ecosystems \cite{Boers2023,Bochow2023}, or in ocean circulation, such as a collapse of the Atlantic meridional overturning circulation \cite{ditlevsen2023warning, van2025collapse}.

With the presence of such tipping points appearing increasingly evident, this has prompted the development of methods for detecting and predicting these state transitions. This has led to a variety of so-called early warning signs that detect approach of a tipping point before it has been crossed \cite{bathiany2025ecosystem, scheffer2009early, Agbedhadji2023}. Among these are physics-based early warning signs that use system process understanding to find system-specific indicators, such as the minimum of the AMOC-induced freshwater transport at the southern boundary of the Atlantic for tipping of the AMOC \cite{van2024physics}, or an increase in the sensitivity of net ecosystem productivity to temperature anomalies that precedes Amazonian rainforest tipping \cite{boulton2013early}. However, the most predominant early warning signs are the generic statistical early warning signs, such as those based on critical slowing down as a system approaches a tipping point. Among these are the methods that measure increases in variance and autocorrelation to detect approach to tipping points \cite{scheffer2009early,Dakos2008}. Despite their widespread use, these methods can sometimes be misinterpreted, as certain events may not exhibit any trend of critical slowing down before a transition occurs \cite{rietkerk2025ambiguity}.

Next to these issues, in general it is also not clear for what kind of transitions these early warning signs do signal. Often, it is assumed that a bifurcation will induce a significant state transition. However, this need not be the case; for example, a Turing bifurcation might lead to the emergence of spatial patterns with limited effects on system functioning -- hereby perhaps even evading tipping altogether \cite{Rietkerk2021, siteur2014beyond, banerjee2023rethinking}. Hence, current statistical early warning signs do not distinguish between proper full system change `Tipping' bifurcations and less dramatic transitions, such as Turing bifurcations.

\begin{figure}
\centering
    \begin{subfigure}[b]{0.49\textwidth}
        \includegraphics[width=\textwidth]{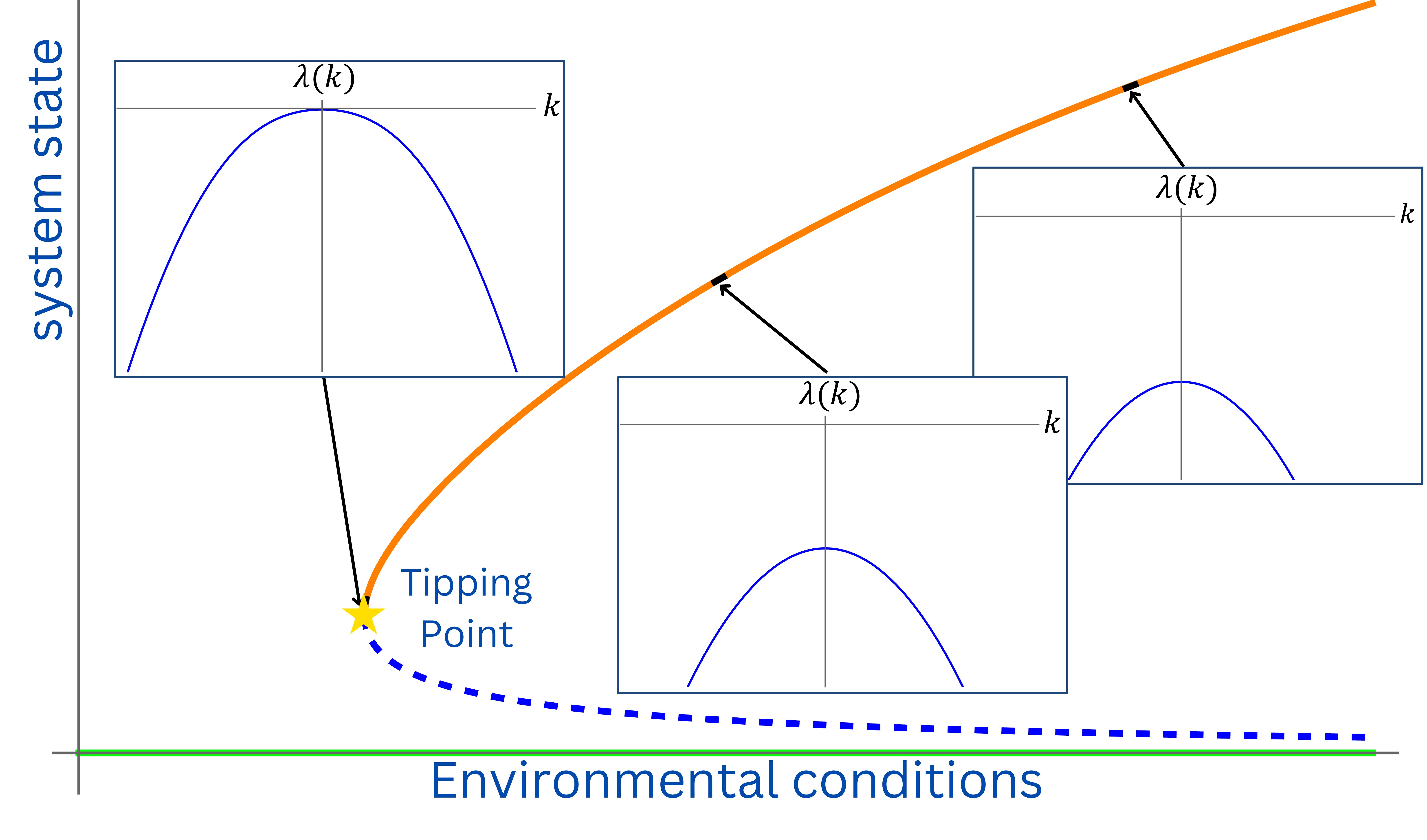}
        \caption{Destabilization via a `Tipping' bifurcation}
        \label{Fig:Tipping_point_Detection}
    \end{subfigure}
    \hfill
    \begin{subfigure}[b]{0.49\textwidth}
        \includegraphics[width=\textwidth]{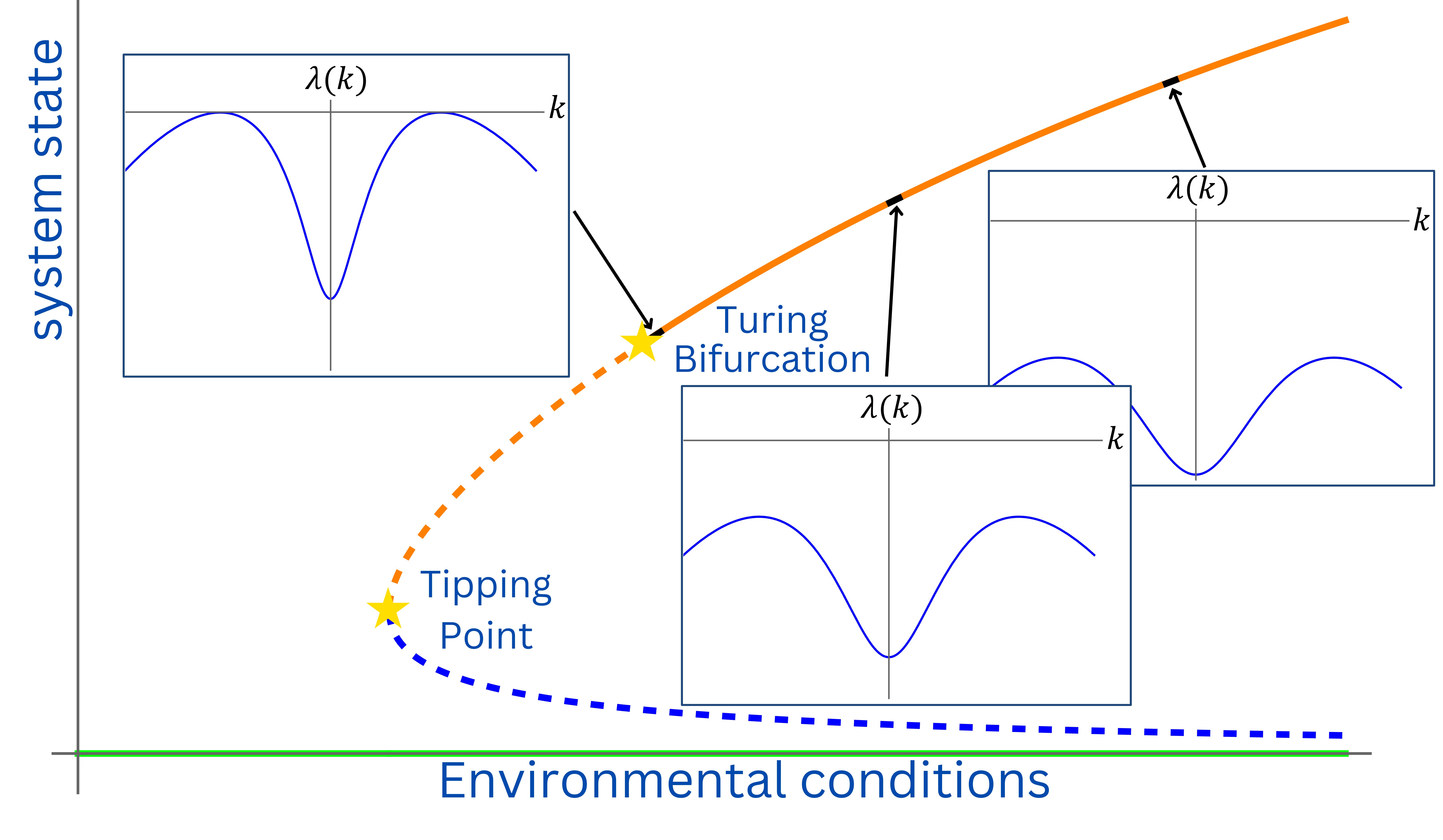}
        \caption{Destabilization via a Turing bifurcation}
        \label{Fig:Turing_bifurcation_Detection}
    \end{subfigure}
    \caption{Illustration of two bifurcation diagrams with distinct destabilizing bifurcations, which are identified using the method introduced in this paper based on the estimation of dispersion relations. The lines indicate the spatially homogeneous steady states, with solid lines indicating stable and dashed lines unstable states. The left shows a situation in which the orange state destabilizes via a saddle-node `Tipping' bifurcation; the right shows a situation in which the it destabilizes via a spatial pattern forming Turing bifurcation. The insets show dispersion relations, relating spatial Fourier eigenmodes, characterized by their wavenumber $k$, to eigenvalues $\lambda(k)$, at various levels with varying closeness to bifurcation. It can be seen that at the bifurcation - and well before it - the dispersion relations are qualitatively different: saddle-node bifurcations are signaled by a peak in the dispersion relation at $k = 0$, whereas Turing bifurcations have peaks for $k \neq 0$. The method introduced in this paper makes use of this fact by providing estimates of these dispersion relations, and thus of the most unstable spatial eigenmodes $k_*$, to determine whether a spatially homogeneous (e.g., left - saddle-node bifurcation) or spatially heterogeneous (e.g., right - Turing bifurcation) destabilization is imminent. Figures are made with the modified Klausmeier model in Eq~\eqref{eq:Extended_klausmeier_2} with parameter values $m = 0.5$, $h = 0.1$, $\delta = 0.5$ (left) or $\delta = 0.01$ (right) and varying $p$ along horizontal axis; system state is represented by the variable $v$.}
    \label{Fig:Early_Warning_Overview}
\end{figure}

In this paper, we introduce a new early warning system that is capable of making this distinction, by providing more information about the kind of bifurcation that is approached (see Figure~\ref{Fig:Early_Warning_Overview}). The method proposed in this paper evaluates the stability of a spatially homogeneous state via an estimation of the associated dispersion relation. This dispersion relation indicates the stability of the system state against different spatially heterogeneous perturbations, specified using spatial (Fourier) modes $k$ and associated mode-dependent eigenvalues $\lambda(k)$. Information on the kind of bifurcation is gained by inferring the critical perturbation that first destabilises the state, i.e. the $k_c$ such that $\mbox{Re}\left(\lambda(k_c)\right) = 0$. If $k_c = 0$, a spatially homogeneous bifurcation occurs (e.g. a saddle-node `Tipping' bifurcation); if $k_c \neq 0$, a spatially heterogeneous bifurcation occurs (e.g. a pattern forming Turing bifurcation). To obtain these central estimates for the dispersion relation, the method fits noisy spatio-temporal data from before a transition to a linear partial differential equation, of which the dispersion relation can be retrieved analytically. 

Here, we investigate the capabilities of this method with synthetic data of a modified-Klausmeier model \cite{bastiaansen2019stable, klausmeier1999regular}, a reaction-diffusion system used to describe the interplay between water and vegetation in drylands. Depending on parameter choices, this model can showcase both full system tipping, organized via a saddle-node bifurcation, and spatial patterning, organized via a Turing bifurcation (see Figure~\ref{Fig:Early_Warning_Overview}). We test how well the dispersion relations can be estimated in various parameter settings, including varying distance to the bifurcation, and for various noise properties. Further, we determine the data needs for this method to work optimal, in terms of requirements on resolution and time span of the input spatio-temporal data.

The rest of this paper is structured as follows. In Section \ref{sec:theory}, we first provide an overview of the theory relating to linear stability of homogeneous states of partial differential equations and dispersion relations. Then, we explain the proposed new early warning system in detail. In Section \ref{sec:Experiment_setup}, we give the details of the numerical experiments performed. In Section \ref{sec:Results}, the results of these experiments are given and discussed. Finally, we end with a brief discussion in Section \ref{sec:discussion}.

\section{Theory} \label{sec:theory}

\subsection{Stability of uniform steady states in spatial systems}\label{sec:theory_theory}

Central to our approach is knowledge about the stability of homogeneous steady states against homogeneous perturbations (signaling e.g. saddle-node bifurcations) and against heterogeneous perturbations (signaling Turing bifurcations). This has been studied mathematically in a wide variety of reaction-diffusion models~\cite{turing1952chemical, siteur2014beyond, doelman2018pattern}. For clarity of presentation in this section we restrict ourselves to the simplest setting of a two-component reaction-diffusion equation, but refer the interested reader to e.g. \cite{doelman2018pattern, siero2020resolving} for extensions to more complicated models.

We consider two-component reaction-diffusion equations with one spatial dimension of the form
\begin{align}\label{eq:theory_genequation}
\begin{split}
    \partial_t u & = \partial_{x}^2 u + f(u,v);\\
    \partial_t v & = \delta \partial_x^2 v + g(u,v),
\end{split}
\end{align}
where $f$ and $g$ represent the reaction terms (i.e. local dynamics) of the model. Let $(u_*,v_*)$ denote a homogeneous steady state satisfying
\begin{align}
    f(u_*,v_*) &= 0, \\
    g(u_*,v_*) &= 0.
\end{align}
We study the stability of this steady state, $(u,v)(x,t) = (u_*,v_*)$, by inspecting the growth rate $\lambda(k)$ of perturbations with wavenumber $k$. To do so, we substitute
\begin{equation}
    \begin{pmatrix}
        u(x,t) \\ v(x,t)
    \end{pmatrix}
    =
    \begin{pmatrix}
        u_* \\ v_*
    \end{pmatrix}
    + e^{\lambda(k) t} e^{i k x}
    \begin{pmatrix}
        \overline{u} \\ \overline{v}
    \end{pmatrix}
\end{equation}
into \eqref{eq:theory_genequation} and linearise the resulting equation. That yields the linear eigenvalue problem
\begin{equation}
    \lambda(k) \begin{pmatrix} \overline{u} \\ \overline{v} \end{pmatrix}
    =
    \begin{pmatrix}
        \frac{\partial f}{\partial u}(u_*,v_*) - k^2 & \frac{\partial f}{\partial v}(u_*,v_*) \\ \frac{\partial g}{\partial u}(u_*,v_*) & \frac{\partial g}{\partial v}(u_*,v_*) - \delta k^2
    \end{pmatrix}
    \begin{pmatrix} \overline{u} \\ \overline{v} \end{pmatrix}
    =: 
    \begin{pmatrix}
        a - k^2 & b \\ c & d - \delta k^2
    \end{pmatrix}
    \begin{pmatrix} \overline{u} \\ \overline{v} \end{pmatrix}
    .\label{eq:theory_lineareq}
\end{equation}
Hence, we obtain two growth rates $\lambda_{1,2}(k)$ as the eigenvalues of the $k$-dependent matrix in \eqref{eq:theory_lineareq}, i.e., as solutions to the $k$-dependent characteristic equation
\begin{equation}\label{eq:Theory_char_eq}
    0 = \lambda^2 + \lambda \left( -a - d + (1+\delta)k^2 \right) + (a-k^2)(d-\delta k^2) - bc.
\end{equation}
These so-called dispersion relations $\lambda_{1,2}(k)$ indicate whether a perturbation with wavenumber $k$ grows or shrinks, i.e. if $\mbox{Re}\left(\lambda_{1,2}(k)\right) < 0$, then perturbations with wavenumber $k$ shrink, and hence the steady state $(u_*,v_*)$ is stable against such perturbations. Thus, if $\mbox{Re}\left( \lambda_{1,2}(k)\right) < 0$ for all wavenumbers $k$, then the steady state $(u_*,v_*)$ is stable against all perturbations. Similarly to ordinary differential equations, there is a bifurcation when (the real part of) one of the eigenvalues changes sign, i.e. when the dispersion relation indicates growth of perturbations of some wavenumber, i.e. when there is a critical wavenumber $k_c$ such that $\mbox{Re}\left( \lambda_1(k_c)\right) = 0$ or $\mbox{Re}\left( \lambda_2(k_c)\right) = 0$. If $k_c = 0$ the bifurcation is due to homogeneous perturbations and indicates, e.g., a saddle-node bifurcation; if $k_c \neq 0$ the bifurcation is due to heterogeneous perturbations and indicates a Turing bifurcation \cite{turing1952chemical}.

As long as the steady state $(u_*,v_*)$ is still stable, the dispersion relation can provide insight into the kind of bifurcation that is approached. In this paper we do so by tracking the most dominant perturbation wavenumber $k_*$, i.e. the wavenumber that corresponds to the slowest decaying regular perturbations. For the system~\eqref{eq:theory_genequation}, it can be deduced that either $k_* = 0$, suggesting approach of a saddle-node bifurcation, or $k_*^2 = \frac{a\delta + d}{2\delta}-\frac{1+\delta}{2\delta}\lambda(k_*)$, suggesting approach of a Turing bifurcation \cite{doelman2018pattern}.

\subsection{Methodology: inferring spatial stability from data}\label{sec:methodology}

\begin{figure}[h!]
\includegraphics[width=\textwidth]{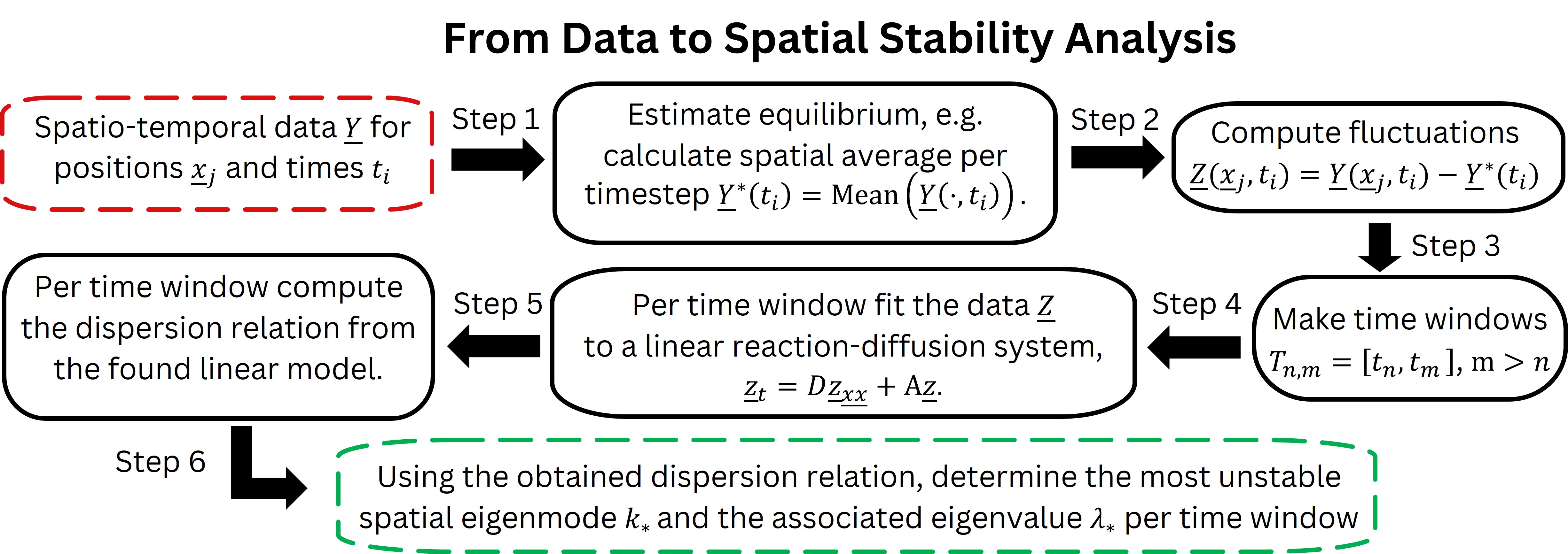}
\caption{Schematic for the proposed method. Spatio-temporal data on fluctuations before a transition is fitted to a linear partial differential equation, whose dispersion relation is computed, and then used to track the most unstable spatial eigenmode $k_*$ and associated eigenvalue $\lambda_* = \lambda(k_*)$.}
\label{Fig:Schematic_Method}
\end{figure}

For the method in this paper we combine linear stability analysis, regression techniques and the study of dispersion relations. An overview is given in the schematic in Figure \ref{Fig:Schematic_Method}. Below, we explain the steps of the method one by one in detail. The general idea is that we use spatio-temporal data from before a transition to determine the fluctuations around an estimated equilibrium state. These are then used to fit to a linear partial differential equation, whose dispersion relations are determined. The implementation for this method as used in this paper is available on \url{https://github.com/JustPaul99/Stability_Analysis_RD}.

As data input, we assume measurements on $\alpha$ variables at times $\{t_1, \ldots, t_{\beta}\}$ and spatial locations $\{ \underline{x}_1, \ldots, \underline{x}_{\gamma}\}$. We denote the full data set by $\underline{Y}$, and a measurement on time $t_i$ and location $\underline{x}_j$ by $\underline{Y}(\underline{x}_j,t_i)$.

\subsection*{Step 1: Computation of the equilibrium}
We approximate the homogeneous equilibrium states as the spatial average of the data per time step, i.e. the equilibrium at time step $t_i$ is estimated as
\begin{align}
    \underline{Y}^*(t_i) := \text{mean}\left(\underline{Y}(\cdot,t_i)\right).
\end{align}
This should be accurate as long as the solution tracks the equilibrium. Other methods to calculate the equilibrium could be chosen such as using data from multiple timesteps (which might be beneficial for example if one does not have many spatial points). 

\subsection*{Step 2: Compute perturbations from equilibrium}
The fluctuations around the estimated equilibrium are calculated per time step, i.e.,
\begin{align}
    \underline{Z}(\underline{x}_j,t_i) :=\underline{Y}(\underline{x}_j,t_i)-\underline{Y}^*(t_i).
\end{align}

\subsection*{Step 3: Divide into time windows}

To analyse the progression of (the stability of) the system, we split the data in time windows, in which we will assume that the dynamics stay the same. These small time windows should be chosen not too big (such that the dynamics stay similar), but also not too small (such that enough data is available). Thus, if the dynamics of a system is changing very slowly, a longer time frame could be effective, and vice versa. For notational parsimony, we omit explicit notation for data in a certain time window.

\subsection*{Step 4: Fit to  linear reaction-diffusion equation}

In a time window, the data $\underline{Z}$ on the fluctuations is expected to behave according to a linear partial differential equation.
\begin{equation}
    \underline{z}_t = A \underline{z} + D \underline{z}_{\underline{x}\underline{x}}, \label{eq:linearPDE}
\end{equation}
where $D := \mbox{diag}(\delta_1,\ldots,\delta_\alpha)$ and
\begin{equation}
    A := \begin{Bmatrix}
a_{1,1} & \cdots & a_{1,\alpha}\\
\vdots&\ddots&\vdots \\
a_{\alpha,1} & \cdots & a_{\alpha,\alpha}
\end{Bmatrix}.
\end{equation}
Now, the goal is to use the data $\underline{Z}$ to find the best estimates for all these parameters. To do so, we rewrite \eqref{eq:linearPDE} to
\begin{equation}
    \underline{z}_t = \Theta \underline{F} := \begin{Bmatrix} A & D \end{Bmatrix} \begin{Bmatrix} \underline{z} \\ \underline{z}_{\underline{x}\underline{x}} \end{Bmatrix}.
\end{equation}
Numerical approximations of the time and spatial derivatives are computed from the data $\underline{Z}$, for example, through finite difference schemes. We denote the thus obtained (estimated) time derivative by $\underline{Z}_t$ and the (estimated) second spatial derivative by $\underline{Z}_{\underline{x}\underline{x}}$. Now, parameter estimation for $\Theta$ can be done, for instance, using least squares fitting, i.e., the optimal $\Theta_*$ is found as
\begin{equation}
    \Theta_* := \mathop{\mathrm{argmin}}_{\Theta} \sum_{i,j} \left\| \underline{Z}_t(\underline{x}_j,t_i) - \Theta \begin{Bmatrix} \underline{Z}(\underline{x}_j,t_i) \\ \underline{Z}_{\underline{x}\underline{x}}(\underline{x}_j,t_i)\end{Bmatrix} \right\|
\end{equation}

\subsection*{Step 5: Compute dispersion relations}
The linear model \eqref{eq:linearPDE} with parameters $\Theta^*$ is analysed using the theory in section \ref{sec:theory_theory}. In particular, using the Fourier mode analysis, i.e. setting $\underline{z}(\underline{x},t) = e^{\lambda(|k|)t}e^{i \underline{k} \cdot \underline{x}} \underline{w}$, dispersion relations $\lambda_m(|\underline{k}|)$ can be computed as the eigenvalues of the matrices
\begin{equation}
    A - |\underline{k}|^2 D = \begin{Bmatrix}
a_{1,1} - \delta_1 |\underline{k}|^2 & a_{1,2} & \cdots & a_{1,\alpha}\\
a_{2,1} & a_{2,2} - \delta_2 |\underline{k}|^2 & \ddots & \vdots \\
\vdots&\ddots&\ddots & a_{\alpha-1,\alpha} \\
a_{\alpha,1} & \cdots & a_{\alpha,\alpha-1} & a_{\alpha,\alpha} - \delta_\alpha |\underline{k}|^2
\end{Bmatrix}.
\end{equation}

\subsection*{Step 6: Tracking most unstable mode}
From the dispersion relations $\{ \lambda_m(|\underline{k}|) \}$, we determine the most unstable wavenumber $k_*$ and its associated (real part of) eigenvalue $\lambda_*$ as
\begin{align}
    k_* &:= \mathop{\mathrm{argmax}}_{|\underline{k}|}\ \max_m\ \operatorname{Re} \left\{ \lambda_m(|\underline{k}|) \right\};\\
    \lambda_* & := \max_{|\underline{k}|}\ \max_m\ \operatorname{Re} \left\{ \lambda_m(|\underline{k}|)\right\}\ = \max_m\ \operatorname{Re}\left\{\lambda_m(k_*)\right\}.
\end{align}
The value of $\lambda_*$ indicates closeness to a bifurcation, as it approaches zero at a bifurcation. The value of $k_*$ is an indicator for the type of bifurcation: if $k_* \approx 0$, it signals a homogeneous bifurcation typically associated with tipping (e.g. a fold bifurcation); if $k_* > 0$, it signals a heterogeneous Turing bifurcation.

\section{Setup of numerical experiments}\label{sec:Experiment_setup}

In this section we will introduce the numerical experiments to test the method from section \ref{sec:methodology}. For these, we have generated synthetic test data using an extended-Klausmeier model for dryland vegetation \cite{bastiaansen2019stable, klausmeier1999regular}, which we introduce in section \ref{sec:method_model}. This model works well for our experiments as it can showcase destabilisations via saddle-node bifurcations and via Turing bifurcations depending on parameter combinations. Subsequently, in subsection~\ref{sec:method_experiments}, we detail the specific numerical experiments that we performed to test the validity and limitations of the method.

\subsection{Test model: an extended Klausmeier model}\label{sec:method_model}
We use the (non-dimensionalised) modified extended Klausmeier model introduced in \cite{bastiaansen2019stable}. This models dryland ecosystems by the interplay between vegetation ($v$) and water ($u$). Here, we use this model on a spatial $1$D domain. The model is given by
\begin{equation}\label{eq:Extended_klausmeier_2}
    \begin{split}
        \d u&=\left(p-u-uv^2+u_{xx}\right)\d t+A\d F^{(1)}_t,\\
        \d v&=\left(uv^2(1-hv)-mv+\delta v_{xx}\right)\d t+A\d F^{(2)}_t.
    \end{split}
\end{equation}
Here, the reaction terms describe the change in water due to rainfall ($+p$), evaporation of water ($-u$), and uptake by plants ($-uv^2$). The change of vegetation is described by mortality ($-mv$) and growth with a carrying capacity ($+uv^2(1-hv)$). Movement of water is modeled as diffusion ($u_{xx}$), and similarly for the dispersal of vegetation ($\delta v_{xx}$). Here, $\delta$ represents the difference in diffusion constants between the two processes (typically $\delta$ is small as vegetation dispersal is slower). Finally, additive spatio-temporal noise $F=(F^{(1)},F^{(2)})^T$ is added, with noise strength $A$. We take $F^{(1)}$ and $F^{(2)}$ as uncorrelated Gaussian processes that are white in time, and either coloured or white in space (depending on the numerical experiment). Here the noise is coloured by applying a normalized squared exponential filter $\exp{(-x^2/l_c^2)}$, with correlation length $l_c$ and position $x$, to white noise\cite{Higdon1998}.

For an analysis of the stability of the homogeneous steady states of the determinstic part of \eqref{eq:Extended_klausmeier_2}, we refer the reader to Appendix~\ref{app:analysis_klausmeier} (and \cite{bastiaansen2019stable}). In short, model \eqref{eq:Extended_klausmeier_2} admits a homogeneous vegetated state $(u_*,v_*)$ for $p > 2m \left( h + \sqrt{1+h^2}\right)$, where
\begin{align}
\begin{split}
    u_* & = m \left( \frac{p}{m} - \frac{v_*}{1-hv_*} \right), \\
    v_* & = \frac{ \frac{p}{m} + \sqrt{ \left(\frac{p}{m}\right)^2 - 4 \left(1 + \frac{p}{m}h\right)}}{2 \left(1+\frac{p}{m}h\right)}.
\end{split}
\end{align}
At $p = p_\textrm{SN} := 2m \left(h + \sqrt{1+h^2}\right)$, there is a saddle-node bifurcation (i.e. tipping point). For $p = p_\textrm{T}$, the Turing bifurcation point, there is no easy closed-form expression, but it can be derived numerically by combining conditions \ref{Eq:Turing_condition_1} and \ref{Eq:Turing_condition_2}. The homogeneous vegetated state is stable for $p > \max\{p_\textrm{SN},p_\textrm{T}\}$ if $p_T$ exists and for $p>p_{SN}$ if $p_T$ does not exist. This means that as $p$ decreases a destabilization occurs either via the aforementioned saddle-node bifurcation or Turing bifurcation.

\subsection{Numerical experiments}\label{sec:method_experiments}
In this study, we take $m = 0.5$ and $h = 0.1$ fixed. To study multiple bifurcation types we typically take either $\delta = 0.5$ or $\delta=0.01$, and let $p$ decrease (i.e. use that as the bifurcation parameter). When $\delta=0.5$, destabilization occurs via a saddle-node bifurcation when $p$ reaches the critical $p_{SN}=1.10499$; when $\delta=0.01$, destabilization occurs via a Turing bifurcation at $p=p_T=1.63398$, which suggests tipping evasion via pattern formation. Bifurcation diagrams for both situations are given in Figure ~\ref{Fig:Early_Warning_Overview}. From hereon, we will refer to parameter settings with $\delta = 0.5$ as the saddle-node case and to parameter settings with $\delta = 0.01$ as the Turing case to remind the reader of the setting studied.

Numerical simulations to generate data were performed as follows. The model \eqref{eq:Extended_klausmeier_2} was spatially discretized using a central difference scheme, and the resulting stochastic differential equation was numerically integrated using the Euler–Maruyama method. We considered two forms of noise: white noise and the aforementioned spatially correlated noise. For the discretizations, we chose a time step of $\text{d}t = 10^{-4}$ and specified the spatial grid size for each of the upcoming experiments in Table~\ref{tab:experiments_final}. We used a finite domain of size $L = 40$ with no-flux boundary conditions.

For this paper, we have designed several numerical experiments to test the applicability and limitations of the method introduced in Section~\ref{sec:methodology}. Specific numerical settings are summarized in Table~\ref{tab:experiments_final}. First, in experiment 1, we investigate how robust the method is against different type of noise, by varying the noise strength and the spatial correlation $l_c$ of the noise. Second, in experiment 2 and 3 we investigate the data requirements for the method by varying the time length of the data, and by varying both the temporal and spatial sampling, which will indicate the kind of resolution that is required for successful employment of the method. Third, in experiment 4, we investigate the robustness of the method in different parameter settings, that vary in closeness to the actual bifurcation (with lower $p$ being closer to the bifurcation), and by also considering the case $\delta = 0.1$ in which Turing and saddle-node bifurcations lie close together (see Figure \ref{fig:Early_warning_shift} for a bifurcation diagram of this situation). Finally, in experiment 5, we investigate a case in which $p$ varies with time to see how the method deals with non-autonomous forcing. For each experimental setting, we construct an ensemble of 100 datasets from model~\eqref{eq:Extended_klausmeier_2}, each with a distinct noise realization. The same set of 100 white noise seeds is used across all different settings, which ensures that the realizations are comparable. We applied the method to each data set separately, and we report on the statistics of the outcomes. 

\begin{table}[h!]
\caption{Summary of numerical settings varied between the numerical experiments. The first row corresponds to baseline values, and subsequent rows to the values for the various numerical experiments. Changes from baseline settings are highlighted in blue ($x$-axes in figures in the results section below) and red ($y$-axes in the figures results section below). Columns denote different numerical settings: $p$ and $\delta$ correspond to parameter values for the model \eqref{eq:Extended_klausmeier_2}; noise strength $A$ and correlation length $l_c$ refer to the Noise properties for the model \eqref{eq:Extended_klausmeier_2}; observation time indicates the length of the time series that went into the method; sample size refers to the resolution of the sampling, in either temporal or spatial component of the data, with 100\% indicating all simulated synthetic data was used; $dx$ refers to the chosen spatial discretizations for the numerical simulation (which was varied in experiment 3 only to better show the effects of spatial sampling). Other settings were kept fixed; specifically, model parameters $m=0.5$ and $h=0.1$, numerical time discretization $dt=10^{-4}$, and spatial domain size $L=40$.}
\label{tab:experiments_final}
\centering

\definecolor{changedcolorx}{HTML}{CCEAFB}
\definecolor{changedcolory}{HTML}{FFD6B3}
\newcolumntype{M}[1]{>{\centering\arraybackslash}m{#1}}

\begin{tabular}{M{2.4cm}|M{1.0cm}M{1.0cm}M{1.2cm}M{1.5cm}M{1.1cm}M{1.3cm}M{1.3cm}M{1cm}}
    \toprule
    \textbf{Experiments} & \boldmath$p$ & \boldmath$\delta$ & \textbf{Noise strength }\boldmath$A$ & \textbf{Corre-lation length } \boldmath$l_c$ & \textbf{Obser-vation time} & \textbf{Sample size time} & \textbf{Sample size space} & \boldmath$\textbf{d} x$ \\
    \midrule
    \rowcolor{gray!30} 
    \textbf{Baseline} & 6 & \shortstack{0.01 \\  0.5} & 1 & 0.1 & 1 & 100\% & 100\% & 0.1 \\
    \midrule
    Experiment 1: Noise & 6 & \shortstack{0.01 \\ 0.5} & \cellcolor{changedcolorx}\shortstack{0.01 \\ 1 \\ 10} & \cellcolor{changedcolory}\shortstack{White \\ 0.1 \\ 0.2} & 1 & 100\% & 100\% & 0.1 \\
    \midrule
    Experiment 2: Temporal sampling & 6 & \shortstack{0.01 \\ 0.5} & 1 & 0.1 & \cellcolor{changedcolory}\shortstack{0.2 \\ 1 \\ 5} & \cellcolor{changedcolorx}\shortstack{100\% \\ 50\% \\ 25\%} & 100\% & 0.1 \\
    \midrule
    Experiment 3: Spatial sampling & 6 & \shortstack{0.01 \\ 0.5} & 1 & 0.1 & 1 & 100\% & \cellcolor{changedcolorx}\shortstack{100\% \\ 50\% \\ 25\%} & \shortstack{0.05} \\
    \midrule
    Experiment 4: Parameters & \cellcolor{changedcolory}\shortstack{2 \\ 6 \\ 20} & \cellcolor{changedcolorx}\shortstack{0.01 \\ 0.1 \\ 0.5} & 1 & 0.1 & 1 & 100\% & 100\% & 0.1 \\
    \midrule
    Experiment 5: Decreasing $p$ & \cellcolor{changedcolory} \shortstack{$20-ct$ \\ \\ c=2} & \cellcolor{changedcolorx} \shortstack{0.01 \\ 0.1 \\ 0.5} & 1 & 0.1 & 1 & 100\% & 100\% & 0.1 \\
    \bottomrule
\end{tabular}
\end{table}

\section{Results}\label{sec:Results}
In this section, we present the results on the aforementioned numerical experiments. We present these results in the form of dispersion relations. Specifically, for each experiment, we show the mean dispersion relation computed from the ensemble, the corresponding $5\%$ and $95\%$ percentiles and the true dispersion relation (based on the parameters of the model). In addition, the densities of the dominant modes $k_*$ and the dominant eigenvalue $\lambda_*$ are estimated using Matlab's \texttt{ksdensity} function and plotted on their respective axes. The critical pairs $(k_*,\lambda_*)$ are shown in a scatterplot, together with an orange ellipse $\sigma(k_*,\lambda_*)$ obtained from their covariance, centered at the mean and aligned with the principal axes with radii given by the standard deviations. In these figures, we also include text boxes displaying the mean and standard deviation of the differences between the estimated and true dominant modes and eigenvalues. Specifically, we report statistics on $\Delta k_* := k_* - k_*^{\text{true}}$ and $\Delta \lambda_* := \lambda_* - \lambda_*^{\text{true}}$, where $k_*, \lambda_*$ refer to estimated modes and eigenvalues, and $k_*^{\text{true}}$, $\lambda_*^{\text{true}}$ to true dominant modes and eigenvalues. To denote these, we use a short-hand notation Mean(standard deviation), but here the order of the standard deviation is presented as the last digit in the mean. For example, $-0.06(36)$ means $-0.06 \pm 0.36$ and $52(47)$ means $52 \pm 47$\cite{EPA2004}. Finally, for completeness, in Appendix \ref{Appendix:Supp_Data}, we report on the statistics of the parameters estimated in step 4 of the method for all of the experiments.

\subsection*{Experiment 1: Effect of noise properties}

\begin{figure}[h!]
\includegraphics[width=\textwidth]{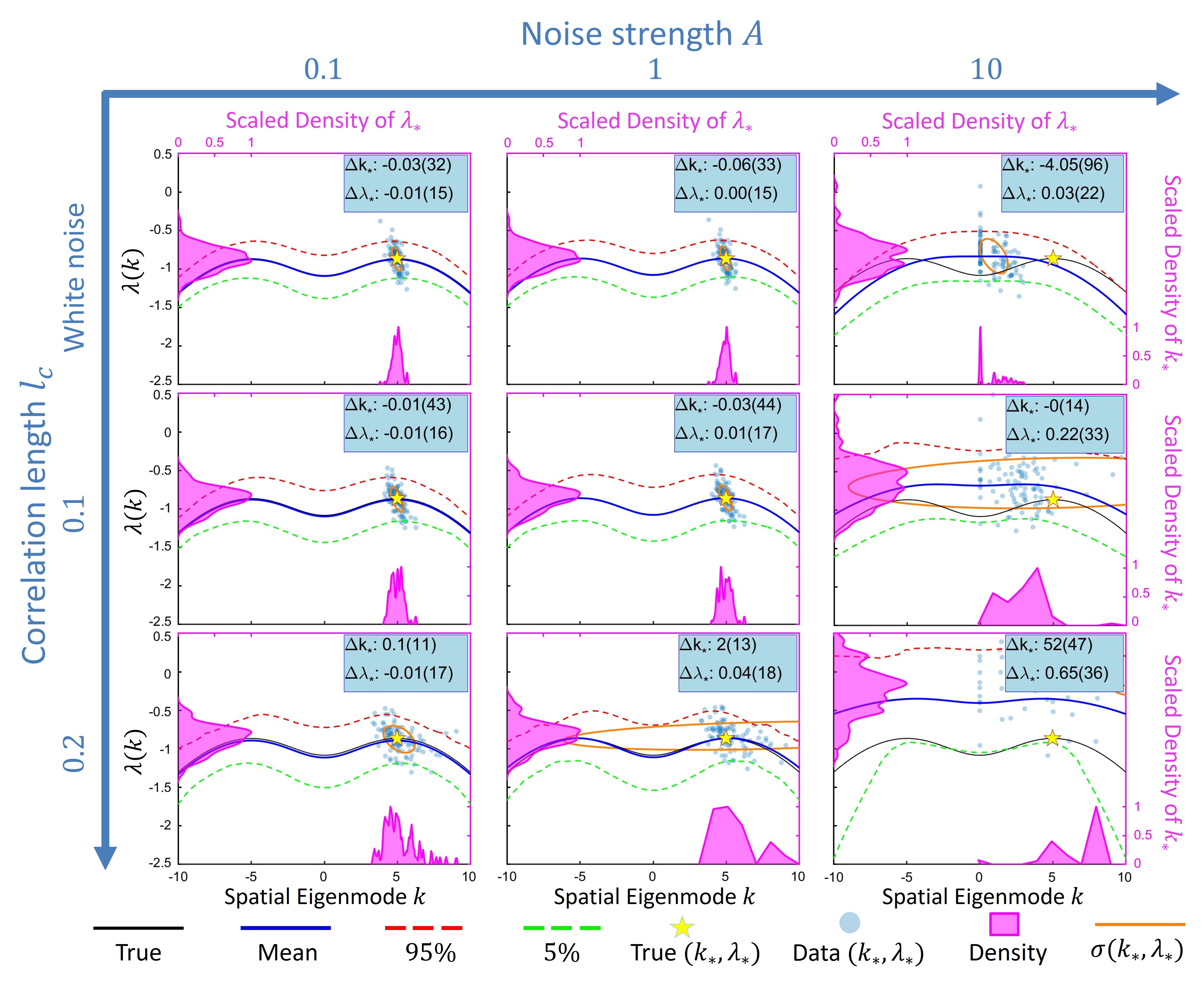}
\caption{Estimated dispersion relations for different noise levels and correlation lengths, with $\delta = 0.5$ (Turing case). The true dispersion relation is shown in black, the ensemble average of estimated dispersion relations in blue, with its 5th and 95th percentiles in red and green dotted lines. The blue dots show the estimated ($k_*$, $\lambda_*$) values from each of the $100$ data sets. The pink areas indicate the scaled marginal densities of these pairs. An orange covariance ellipse, centered at the mean, represents the one-standard-deviation, $\sigma(k_*,\lambda_*)$, spread along the principal directions. The true dominant pair ($k_*^{\text{true}}$, $\lambda_*^{\text{true}}$) is shown with a star. Blue shaded insets report on statistics of error of the method in determining this dominant pair (see main text). The data has been generated using the parameter settings defined in Experiment 1 from table \ref{tab:experiments_final}.}
\label{Fig:Noise_level_correlation_Turing}
\end{figure}

\begin{figure}[h!]
\includegraphics[width=\textwidth]{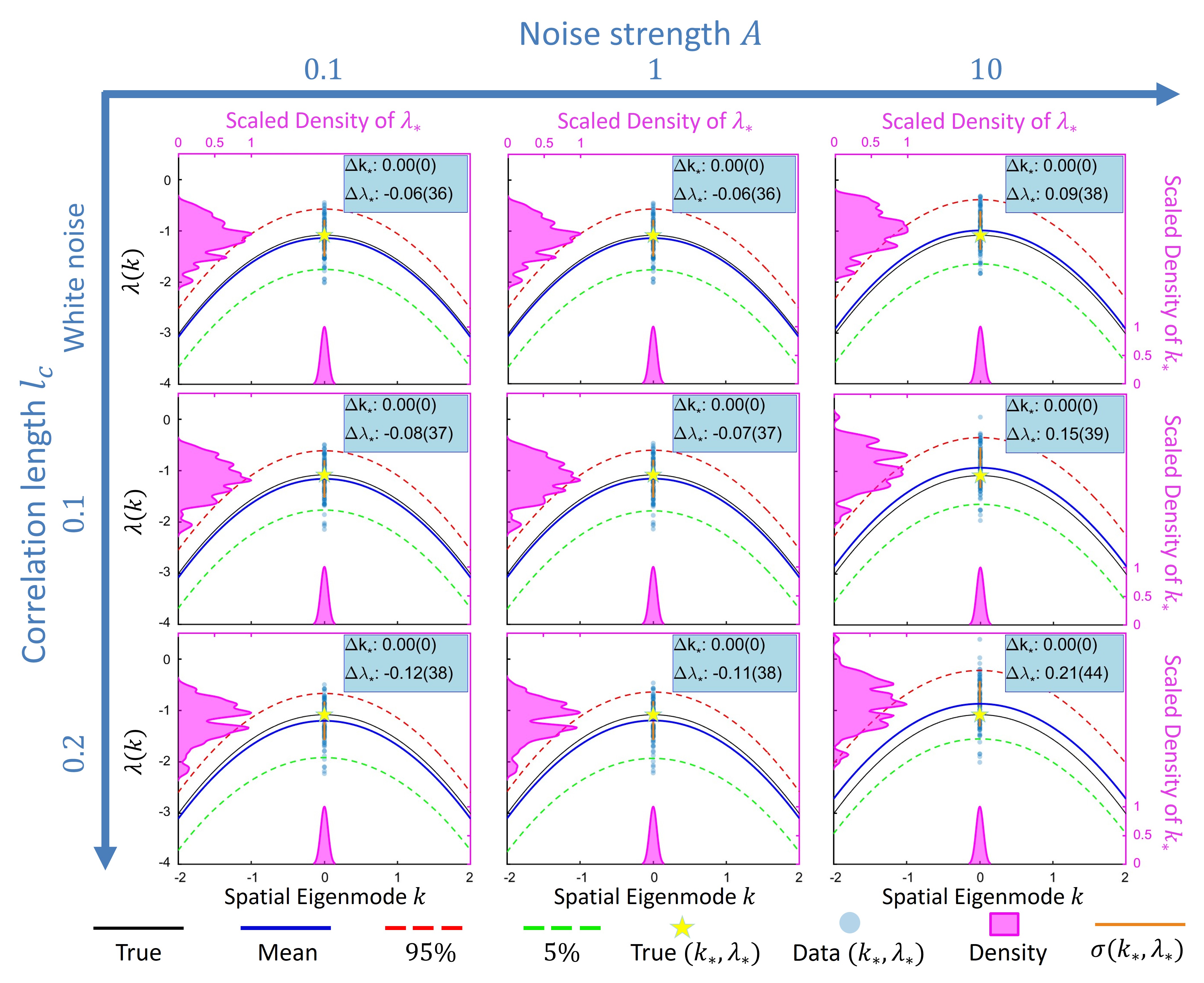}
\caption{Estimated dispersion relations for different noise levels and correlation lengths, with $\delta=0.5$ (saddle-node case). Data is generated using Experiment 1 settings in Table~\ref{tab:experiments_final}. See the caption of Figure \ref{Fig:Noise_level_correlation_Turing} for the details of the depicted lines, areas, circles and insets}
\label{Fig:Noise_level_correlation_Tipping}
\end{figure}

Figures \ref{Fig:Noise_level_correlation_Turing} and \ref{Fig:Noise_level_correlation_Tipping} provide an overview of how the method responds to varying noise levels and correlation lengths. For low noise strengths ($0.1$ and $1$) and short correlation lengths (White noise and $0.1$), the dominant spatial eigenmodes identified by the method show little variation, and the estimated dispersion relations are in line with the true one. However, the method fails under higher noise strength and larger correlation lengths. For example, in Figure \ref{Fig:Noise_level_correlation_Turing} the method fails to recover the correct dispersion relation at a noise strength of $10$. As shown in Figure~\ref{fig:Heatmap_effect_Noise}, high noise levels drive the vegetation below the unstable homogeneous vegetation state. We suspect this indicates that the system is locally being pulled toward alternative attractors. In those cases the method cannot capture all the relevant (nonlinear) dynamics. 

Moreover, we find that the spread in the estimated dominant spatial eigenmode gets bigger as the noise's correlation length increases. We suspect that the interplay of the higher correlation in the noise can make the method unreliable for detecting some noise realizations. For instance, Figure~\ref{fig:correlation_0.2} shows all 100 estimated dispersion relations for the Turing case ($\delta = 0.01$) at noise level $0.1$ and correlation length $l_c = 0.2$, where some realizations exhibit peaks at different wavenumbers.

Finally, we observe that the saddle-node case with $\delta = 0.01$ exhibits greater robustness to varying noise properties compared to the Turing case with $\delta = 0.5$. We suspect this is because retrieval of the dominant homogeneous spatial eigenmode is less sensitive to the spatial structure of noise. In fact, the fitted parameters of the linear model in step 4 of the model become less accurate for higher noise strength (see Table \ref{Table:Noise_level_correlation_Tipping}).

\subsubsection*{Experiment 2 \& 3: Data requirements}
\begin{figure}[h!]
\includegraphics[width=\textwidth]{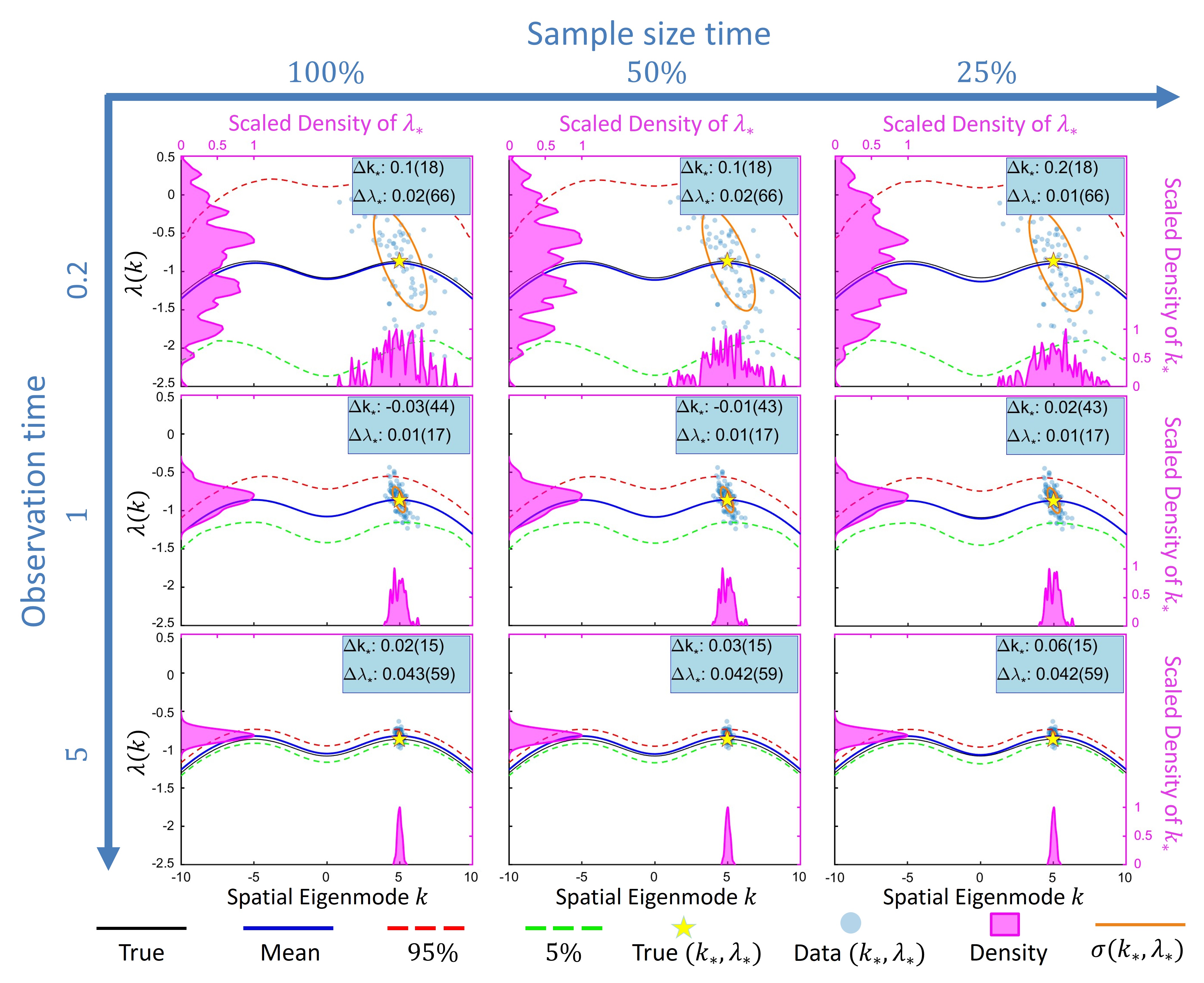}
\caption{Estimated dispersion relations for different observation times and temporal sampling rates for $\delta = 0.01$ (Turing case). Data is generated using Experiment 2 settings in Table~\ref{tab:experiments_final}. See the caption of Figure \ref{Fig:Noise_level_correlation_Turing} for the details of the depicted lines, areas, circles and insets.}
\label{Fig:sampling_Turing}
\end{figure}

\begin{figure}[h!]
\includegraphics[width=\textwidth]{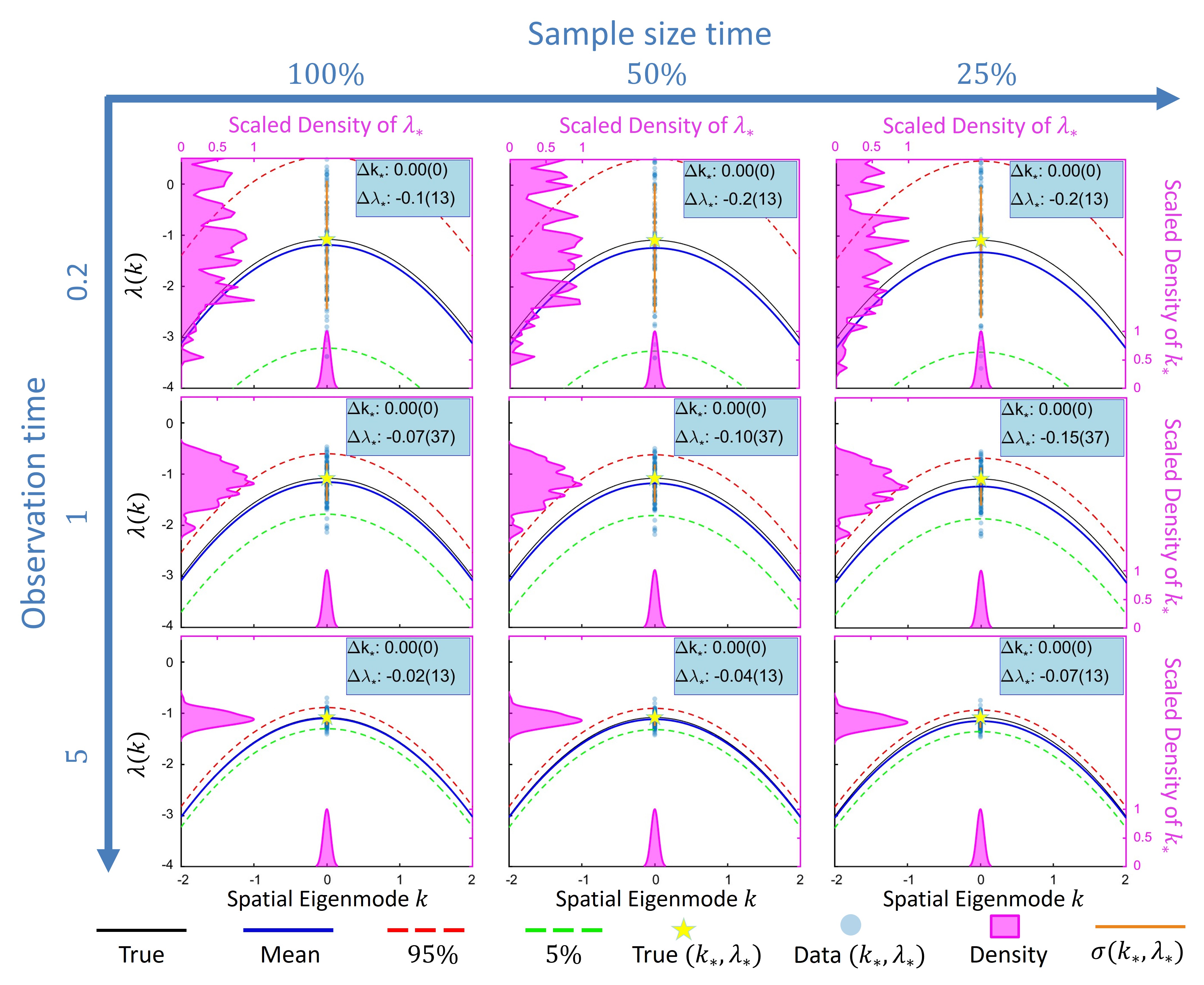}
\caption{Estimated dispersion relations for different observation times and temporal sampling rates for $\delta = 0.5$ (saddle-node case). Data is generated using Experiment 2 settings in Table~\ref{tab:experiments_final}. See the caption of Figure \ref{Fig:Noise_level_correlation_Turing} for the details of the depicted lines, areas, circles and insets.}
\label{Fig:sampling_Tipping}
\end{figure}

Figures~\ref{Fig:sampling_Turing} and~\ref{Fig:sampling_Tipping} show the results of experiment 2 in which we vary the total observation time and the temporal sampling rate for the Turing and saddle-node case respectively. These show that increasing the total observation time, and hence using more data, improves the estimations; the spread in dominant spatial eigenmodes $k_*$ and eigenvalues $\lambda_*$ goes down with longer observation times, and the spread in estimated dispersion relations goes down as well. Further, with respect to temporal sampling, it can be seen by eye that the estimated dispersion relations do not change much when less of the time steps have been used. However, the errors for $k_*$ and $\lambda_*$ do increase a bit (as do the fitted parameters in step 4; see Table \ref{Table:Sampling_Turing}). In general, it can thus be seen that the total observation time, rather than the amount of data, leads to most improvements.

\begin{figure}[h!]
\includegraphics[width=\textwidth]{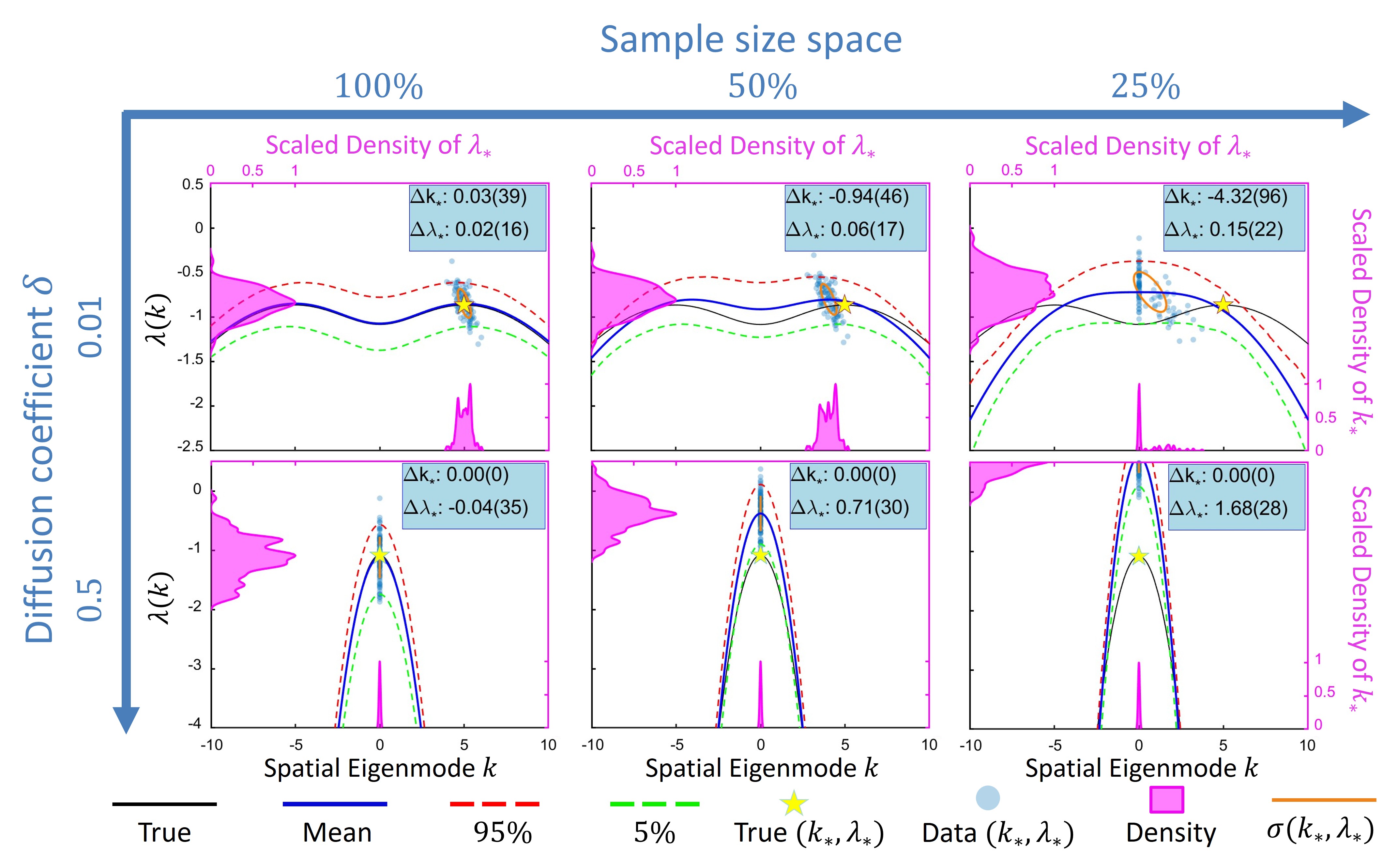}
\caption{Estimated dispersion relations for different spatial sampling rates for $\delta = 0.01$ (Turing case; top row) and $\delta = 0.5$ (saddle-node case; bottom row). Data is generated using Experiment 3 settings in Table~\ref{tab:experiments_final}. See the caption of Figure \ref{Fig:Noise_level_correlation_Turing} for the details of the depicted lines, areas, circles and insets. Note that in the bottom row different scales have been used compared to Figure~\ref{Fig:Noise_level_correlation_Tipping} and \ref{Fig:sampling_Tipping}.}
\label{Fig:spatial_sampling_corr_02}
\end{figure}

Figure \ref{Fig:spatial_sampling_corr_02} shows the results of experiment 3 in which the spatial sampling rate is varied for both the Turing and the saddle-node case. We see that lower spatial sampling rates worsen the estimation -- in particular, the estimation of the dominant spatial eigenmodes $k_*$. This stems from the fact that spatial sampling reduces the accuracy of the estimation of the diffusion process (see Table \ref{Table:spatial_sampling_corr_02}). During testing we observed that the correlation length of the noise does influence how much spatial sampling is acceptable: the longer the noise's correlation length, the more coarse spatial sampling could be to still have meaningful estimations. For example, Figure~\ref{Fig:spatial_sampling_white} and Table~\ref{Table:spatial_sampling_white} illustrate that, under spatial white noise, spatial sampling is effectively not possible.

\subsubsection*{Experiment 4: sensitivity to parameters}

\begin{figure}[h!]
\includegraphics[width=\textwidth]{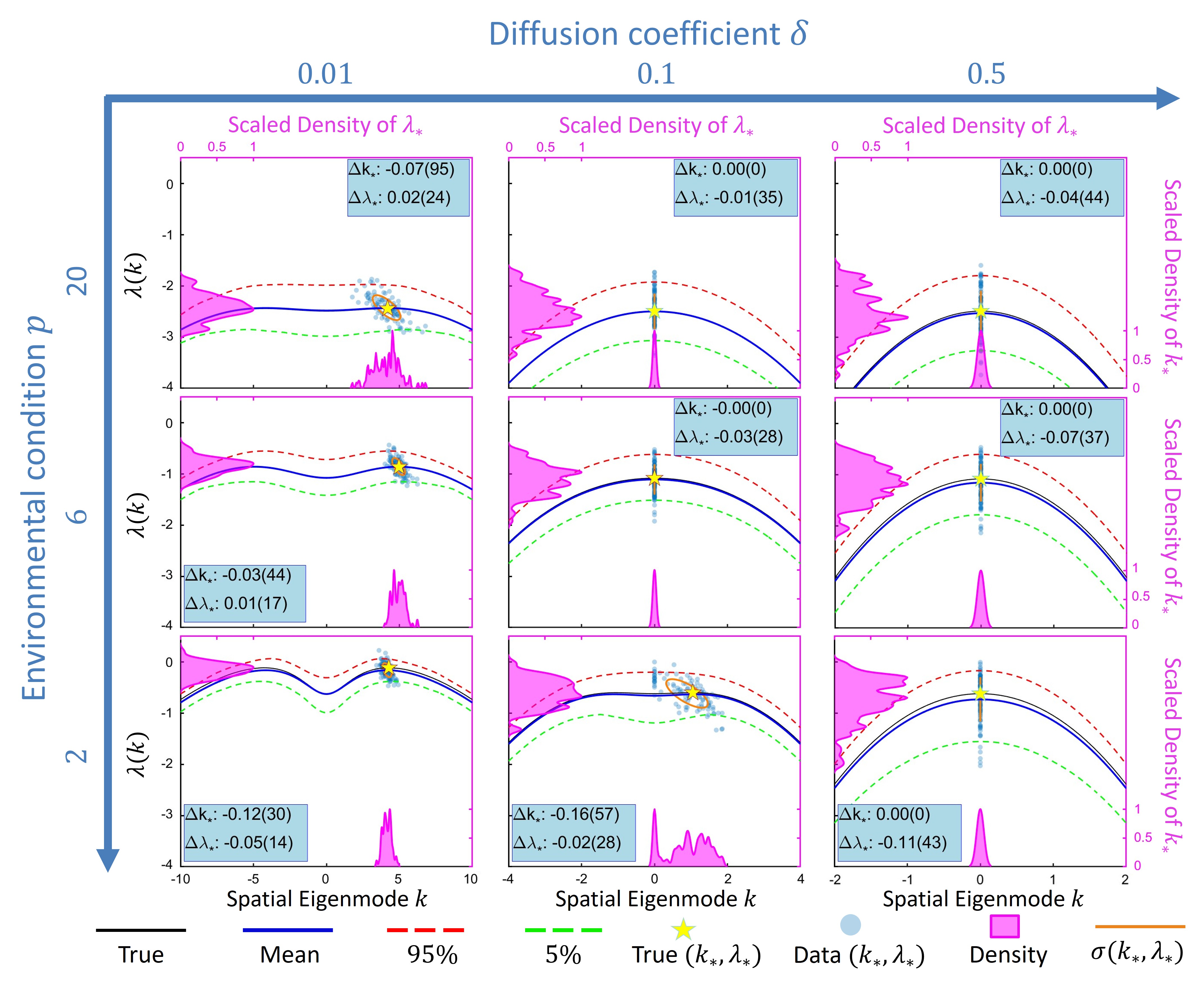}
\caption{Estimated dispersion relations for different parameters $\delta$ and $p$. Data is generated using Experiment 4 settings in Table \ref{tab:experiments_final}. See the caption of Figure \ref{Fig:Noise_level_correlation_Turing} for the details of the depicted lines, areas, circles and insets. Note that the scales in this figure are different compared to the previous figures.}
\label{Fig:parameter_tests}
\end{figure}

Figure~\ref{Fig:parameter_tests} shows the results of Experiment 4, in which the model parameters $p$ and $\delta$ were varied. Lower values of $p$ bring the system closer to a bifurcation, with critical values $p_c$ given by: $p_c=1.63398$ for $\delta=0.01$ (Turing), $p_c=1.11445$ for $\delta=0.1$ (Turing), and $p_c=1.10499$ for $\delta=0.5$ (saddle-node).
These figures show that the method correctly finds the dispersion relations and dominant spatial eigenmodes $k_*$ and eigenvalues $\lambda_*$ in many cases. In particular, the closer the system is to a bifurcation (i.e. lower $p$), the better the prediction of the dominant spatial eigenmode $k_*$ corresponds to the critical spatial eigenmode $k_c$ at the true bifurcation. 

However, the method seems to have difficulties in the specific situation near a shift between a dominant homogeneous spatial eigenmode ($k_*=0$) and a dominant heterogeneous spatial eigenmode ($k_*\neq 0$), as seen for parameters $p=20,\delta=0.01$ and $p=2,\delta=0.1$. In these cases, the dispersion relation is near horizontal, which leads to increased sensitivity resulting in higher errors for estimated dominant eigenmodes $k_*$.

Further, it can be seen in some of the panels that sometimes an eigenvalue with positive real part is unjustly estimated. In these cases, increasing the observation time would lead to more accurate estimations, per experiment 2, removing this issue. However, current results would still accurately indicate the dominant spatial eigenmode $k_*$, and thus the type of bifurcation (except in aforementioned situation when saddle-node and Turing bifurcations happen almost simultaneously).

\subsubsection*{Experiment 5: Time-dependent forcing}
\begin{figure}[h!]
\includegraphics[width=\textwidth]{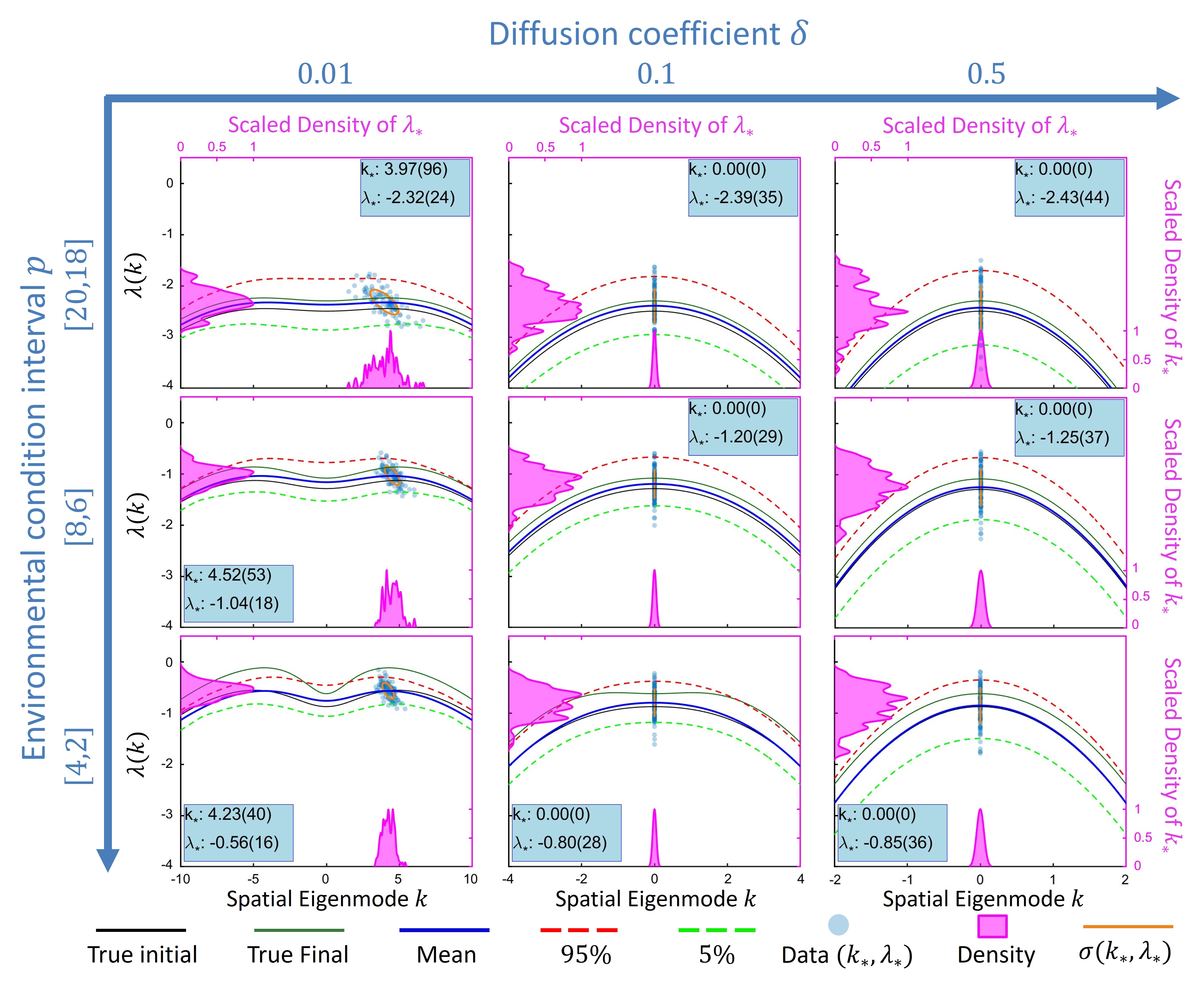}
\caption{Estimated dispersion relations for different parameters $\delta$ an time-varying parameter $p(t) = 20 - 2t$ with measurements between $p=20$ and $p=18$ (top row), between $p = 8$ and $p = 6$ (middle row) and between $p = 4$ and $p = 2$ (bottom row). Data is generated using Experiment 5 settings in Table \ref{tab:experiments_final}. The black lines indicate the true initial dispersion relations (i.e. at the highest $p$ value of the interval), and the darkgrey lines the true final dispersion relations (i.e. at the lowest $p$ value of the interval). In contrast to the previous figures, the blue insets denote the values of $k_*$ and $\lambda_*$. See the caption of Figure~\ref{Fig:Noise_level_correlation_Turing} for the details of the other depicted lines, areas, circles and insets.}
\label{Fig:parameter_time_variatiations}
\end{figure}

In the final experiment, we have introduced a time-varying parameter $p$, where $p$ goes down from $p=20$ to $p=2$ linearly as $p(t)=20-2t$. In Figure \ref{Fig:parameter_time_variatiations} the results are shown when the methodology is used on time windows between $p = 20$ and $p = 18$ (top row), between $p = 8$ and $p = 6$ (middle row), and between $p =4$ and $p=2$ (bottom row). The estimated dispersion relations are similar to those in Figure \ref{Fig:parameter_tests} with constant parameters. Hence, this indicates that the method is capable of handling time-varying parameters.

In general, we find that the estimated dispersion relation lies roughly between the dispersion relations for the initial and final parameter values of the given parameter window. Similarly to the constant parameter case, the dominant spatial eigenmode is estimated well in most situations, except when the dispersion relation is nearly horizontal (e.g., $\delta = 0.01$, $p \in [20,18]$; and $\delta = 0.1$, $p \in [4,2]$). For $\delta = 0.01$, $p \in [4,2]$, the dispersion relation changes substantially between the start ($p = 4$) and end ($p = 2$), which is not fully captured by the method. However, the method correctly distinguishes between a dominant spatial eigenmode with $k_* = 0$ and $k_* \neq 0$. As in the constant parameter case, the closer the system is to crossing a bifurcation (i.e., lower $p$), the more accurately the estimated $k_*$ indicates whether a Turing or saddle-node bifurcation is approached — except when the bifurcations are very close together (e.g., $\delta = 0.1$).

\section{Discussion}\label{sec:discussion}

In this paper, we have introduced a new early warning system designed to distinguish between imminent spatially homogeneous `Tipping' destabilizations and imminent spatially heterogeneous `Turing' destabilizations. For this, spatio-temporal data from fluctuations before the crossing of a bifurcation is fitted to a linear reaction-diffusion system. Its associated dispersion relation, that relates spatial Fouriermodes $k$ to eigenvalues $\lambda(k)$, then yields the most dominant spatial mode $k_*$, which is used to distinguish between tipping $(k_c=0)$ and Turing $(k_c \neq 0)$ destabilizations. 

Using numerical experiments on synthetic data, we have shown the effectiveness of this methodology in estimating dispersion relations and making the distinction between Tipping and Turing bifurcations. In particular, the results indicate the methodology can handle many different noise dynamics and also time-varying forcing. When used on appropriate data, the method seem to fail only in two cases. First, when noise is too dominating (and e.g. noise-induced transitions occur). Second, when the dispersion relation is very flat (e.g. due to switching between Tipping and Turing destabilisation).

These numerical experiments also highlighted the data requirements to effectively use this methodology. Foremost, it has been indicated that it works better with longer observation times; however, the temporal resolution of the data was less important. This indicates that longer, but potentially less frequent, measurements should be used to improve the estimations -- especially for the eigenvalues. At the same time, the spatial resolution of the data is in fact very important, but the specific resolution needs depend on the spatial correlation of the noise. Therefore, when using this method, measurements should be used that have a spatial resolution at least as fine as the expected spatial correlation in the system. So for systems were disturbances influence large areas - such as forest fires - a coarser spatial resolution might still work whereas for systems with primarily small-scale disturbances it will not.

Within this paper, we focused on synthetic data from a model in which there could be a saddle-node `Tipping' bifurcation or a Turing bifurcation. However, this method should also be able to detect destabilizations of spatially homogeneous states organized by other bifurcations, such as pitchfork or transcritical bifurcations. Further, also Hopf or Turing-Hopf bifurcations can be detected by explicitly taking the imaginary part of the eigenvalues into account. Extensions to other systems, with potentially more complicated dispersion relations, also should be possible by e.g. incorporating higher order spatial derivatives or more components in step 4 of the method. Also for spatially heterogeneous states, the method might be useful, although the stability of such states is typically not fully described by dispersion relations and additional analysis would be needed \cite{sandstede2002stability}.

The method developed here is capable of making statements about the linear stability of a system state, and how it might destabilize. However, it cannot make statements about the new state after crossing of any bifurcation. For example, while it has been argued that a Turing bifurcation can lead to small-amplitude spatial patterns and evastion of tipping \cite{Rietkerk2021}, spatially heterogeneous destabilisations can also lead to large-amplitude patterns or even expedite tipping \cite{van2025vegetation, pinto2025spatial}. Further analysis that incorporates the nonlinear feedbacks of a system is needed to make such statements. For example, to differentiate between super-critical Turing (small amplitude patterns emerge) and sub-critical Turing (large amplitude patterns emerge) bifurcations, fits to amplitude equations, that incorporate the nonlinear effects, might be fruitful \cite{doelman2018pattern}.

Finally, at the heart of the early warning system is the fit to a linear model in step 4. In this study, we have used linear regression and numerical approximations of derivatives, but more refined approaches could further improve the results. Specifically, for stochastic spatial systems, progress has been made in retrieving models directly from data \cite{Boninsegna2018,Callaham2021}. Additional improvements could be achieved using alternative regression methods \cite{Brunton_2016,rudy2017data}, bootstrapping \cite{Fasel_2022}, weak formulations of the partial differential equations \cite{Fasel_2022,MESSENGER2021110525}, or by explicitly accounting for the time-dependence of the forcing \cite{Messenger2022,Fasel_2022}.

Understanding the timing and nature of ecosystem and climate subsystem destabilization is a pressing challenge in the face of ongoing climate change. The methodology presented in this paper addresses this challenge by estimating dispersion relations from spatio-temporal data collected before a transition. In this way, both the eigenvalues and the dominant spatial eigenmodes are determined, providing information on the timing of destabilization as well as the type of bifurcation. Hence, this methodology might form the basis of refined statistical early warning system that can signal not only when, but also what happens at a critical transition.

\section*{Data statement}
All codes used in this paper are available on \url{https://github.com/JustPaul99/Stability_Analysis_RD}.

\section*{Declaration of generative AI and AI-assisted technologies in the writing process}
During the preparation of this work the author(s) used ChatGPT, Gemini, and Copilot to improve readability of the text. After using this tool/service, the author(s) reviewed and edited the content as needed and take(s) full responsibility for the content of the published article.

\bibliographystyle{unsrturl}
\bibliography{References}

\newpage

\appendix
\renewcommand{\thefigure}{A.\arabic{figure}}
\renewcommand{\thetable}{A.\arabic{table}}
\renewcommand{\theequation}{A.\arabic{equation}}
\setcounter{figure}{0}
\setcounter{equation}{0}
\setcounter{table}{0}

\section{Analysis of extended Klausmeier model~\eqref{eq:Extended_klausmeier_2}}\label{app:analysis_klausmeier}
In this section, we perform a stability analysis of the spatially homogeneous states of the (deterministic part of the) model \eqref{eq:Extended_klausmeier_2}. Conform \cite{bastiaansen2019stable}, the system has three spatially homogeneous steady states: the no-vegetation solution $(u_0,v_0)=(p,0)$, and two vegetation equilibria $(u_{1,2},v_{1,2})$ given by
\begin{align}
\begin{split}
    u_{1,2}&=\frac{2 hm + p +2h^2p\pm\sqrt{p^2-4hmp-4m^2}}{2+2h^2},\\
    v_{1,2}&=\frac{p\mp\sqrt{p^2-4hmp-4m^2}}{2m+2hp}.
\end{split}
\end{align}
These vegetation states exist for $p>p_{SN}:=2m(h + \sqrt{1+h^2})$ as we can see that for both $u_{1,2}$ and $v_{1,2}$ the square root needs to be positive for the solutions to exist. 

First, we determine their stability against spatially homogeneous equations, i.e. by inspecting the Jacobian \eqref{eq:theory_lineareq} for $k=0$,
\begin{align}
A_{i}(0)=
    \begin{pmatrix}
        -1-v_i^2 & -2v_iu_i\\
        v_i^2(1-hv_i) & - m +2 u_i v_i -3 hu_iv_i^2
    \end{pmatrix}.
\end{align}
The no vegetation state $(u_0,v_0)$ is stable for all $p\geq0$, as $A_0(0)$ has eigenvalues $\lambda_1=-1<0$ and $\lambda_2=-m<0$. For the vegetation states we check the stability conditions $\det(A_{1,2}(0))>0$ and $\Tr(A_{1,2})<0$. After some calculations we find
\begin{align}
    \det(A_{1,2})(0)=\frac{D\mp p\sqrt{D}}{2m+hp},
\end{align}
with $D=-4 m^2-4 h m p+p^2$. In the domain for $p>p_{SN}$ $\det(A_1)$ is positive and $\det(A_2)$ is negative. So $(v_1,w_1)$ is unstable and we check the stability of $(v_2,w_2)$ with the trace. For the trace after some calculations, we find
\begin{align}
    \Tr(A_2)=-1-v_2^2+u_2v_2(1-2hv_2)<0.
\end{align}
This is not a trivial condition to uphold, but following Reduce from Mathematica the condition $0<m\leq2+2h^2$ assures that the trace remains negative. So the vegetation state $(u_2,v_2)$ is stable for $p>p_{SN}$ and $0<m\leq 2+2 h^2$.

For the simulations introduce in section \ref{sec:Experiment_setup}, we use models with $m<2$, and we let $p>p_{SN}$ be our environmental condition that we change to induce tipping or Turing patterns. Therefore, stability on the homogeneous spatial patterns is satisfied.

Second, we perform an analysis on when we should expect tipping behavior or Turing patterns to emerge. We can find the critical eigenmode $k_c$ by differentiating over the $k$-dependent characteristic equation \ref{eq:Theory_char_eq} with respect to $k$ and solving for $k_c$. Here we can use $\frac{\d\lambda}{\d k}|_{k_c,\mu_c}=0$ and $\lambda|_{k_c,\mu_c}=0$, where $\mu_c$ represents the critical parameters, i.e. $p_{T}$ for a fixed $h$, $m$ and $\delta$. This gives $k_c=\sqrt{\frac{{a\delta+d}}{2\delta}}$, where $a$ and $d$ are as defined in \eqref{eq:theory_lineareq}

The critical parameter $p_T$ can then be found by substituting $k_c$ in the characteristic equation \eqref{eq:Theory_char_eq} to find the indirect nontrivial relation,
\begin{align}
    (a\delta-d)^2=-4\delta b c,
\end{align}
where $a,b,c,d$ are as defined in \eqref{eq:theory_lineareq}. Using the above relation one can find $p_T$ from this relation once $h$, $m$ and $\delta$ are fixed. Then from $p_T$, $h$ and $m$ we can calculate the critical wavenumber, $k_c$.

For Turing patterns, we require that $k_c>0$, which implies that $a\delta+d>0$. Together with the stability conditions for the homogeneous spatial vegetation state, $a+d>0$ and $ad-bc<0$, from the trace and determinant respectively, we find two minimal conditions for Turing patterns within a stable vegetation state:
\begin{align}
    a<0, \quad d>0,\quad bc<0, \quad 0<\delta<1 \quad \textbf{OR} \quad 
    a>0, \quad d<0, \quad bc<0 \quad \delta>1.
\end{align}

For the Extended Klausmeier model we have that $a_{KL}=-1-v_2^2<0$, $b_{KL}=-2v_2u_2<0$, $c_{KL}=v_2^2(1-hv_2)>0$\footnote{Follows from substitution of the minimal parameter $h_{min}=\frac{-2m+p}{p}$ in $(1-hv_1)|_{h=h_{ex}}=\frac{2m}{p}>0$. } and $d_{KL}=u_2v_2(1-2hv_2)$, the sign of which is indeterminate in general. $d_{KL}$ is therefore the only function we can play with to study different settings for the emergence of saddle-node bifurcations and Turing bifurcations. This means that for $d_{KL}<0$ the spatially homogeneous vegetation state is always stable. However, as we study $\delta<1$ and $a_{KL}<0$, the emergence of Turing patterns or saddle-node bifurcations can only occur when $d_{KL}>0$.  Altogether, the conditions for Turing patterns are
\begin{align}\label{Eq:Turing_condition_1}
d_{KL}|_{p_T}>\delta |a_{KL}||_{p_T}
\end{align}
with $p_T^2-4hmp_T-4m^2>0$, i.e. $p_T>p_{SN}$, where $p_T$ follows from
\begin{align}\label{Eq:Turing_condition_2}
    (\delta a_{KL}-d_{KL})^2=-4\delta b_{KL} c_{KL}.
\end{align}
\newpage

\section{Supplementary Data: fitted model parameters} 
\label{Appendix:Supp_Data}
\renewcommand{\thefigure}{B.\arabic{figure}}
\renewcommand{\thetable}{B.\arabic{table}}
\renewcommand{\theequation}{B.\arabic{equation}}
\setcounter{figure}{0}
\setcounter{equation}{0}
\setcounter{table}{0}

In this section, we report on additional results for the numerical experiments. These includes reporting on the found fitted parameters in step 4 of the method. In this paper, the data of fluctuations $\underline{Z} = \begin{Bmatrix} \bar{u} & \bar{v} \end{Bmatrix}^T$ was fitted to the linear model (per equation~\ref{eq:linearPDE})
\begin{align}
\begin{split}\label{eq:linearPDE_app}
    \partial_t \bar{u} & = \partial_{xx} \bar{u} + a \bar{u} + b \bar{v}\\
    \partial_t \bar{v} & = \delta \partial_{xx} \bar{v} + c \bar{u} + d \bar{v}. 
\end{split}
\end{align}
To denote the estimated parameters of this linear model, we use -- as in the main text -- the short-hand notation Mean(standard deviation), but here the order of the standard deviation is presented as the last digit in the mean. For example, $-0.06(36)$ means $-0.06 \pm 0.36$ and $52(47)$ means $52 \pm 47$\cite{EPA2004}.

\subsection{Additional Material for Numerical Experiment 1}\label{Appendix:Experiment_1}

\begin{table}[H] 
\begin{center}
\caption{Estimated parameters of the linear model \ref{eq:linearPDE_app} to which spatio-temporal data of fluctuations were fitted. This table denotes the values for the same settings as shown in Figure \ref{Fig:Noise_level_correlation_Turing}, i.e. Experiment 1 from Table \ref{tab:experiments_final} with $\delta=0.01$ (Turing case).} 
\label{Table:Noise_level_correlation_Turing}

\begin{tabular}{@{}c@{\hskip 0.5em}c@{}}
\raisebox{-3.5em}{\rotatebox{90}{\textbf{Correlation Length}}} & \shortstack{\textbf{Noise Level} \\
\begin{tabular}{|c|c|c|c|}
\toprule
\textbf{} & 0.1 & 1 & 10 \\
\midrule
\textbf{\rotatebox{90}{White Noise}} &
\shortstack{
$a: -29.79(15)$\\
$b: -2.16(15)$\\
$c: 13.26(15)$\\
$d: -0.09(15)$\\
$\delta: 0.01(15)$
} &
\shortstack{
$a: -29.61(15)$\\
$b: -2.16(15)$\\
$c: 13.06(15)$\\
$d: -0.09(15)$\\
$\delta: 0.01(15)$
}&
\shortstack{
$a: -29.74(20)$\\
$b: -1.20(22)$\\
$c: 7.12(20)$\\
$d: -0.54(22)$\\
$\delta: 0.01(22)$
}\\
\midrule
0.1 & \shortstack{
$a: -29.78(16)$\\
$b: -2.16(16)$\\
$c: 13.25(16)$\\
$d: -0.10(16)$\\
$\delta: 0.01(16)$
} & \shortstack{
$a: -29.38(17)$\\
$b: -2.17(17)$\\
$c: 12.93(17)$\\
$d: -0.08(17)$\\
$\delta: 0.01(17)$
}
 & \shortstack{
$a: -27.87(30)$\\
$b: -0.91(33)$\\
$c: 8.64(30)$\\
$d: -0.40(33)$\\
$\delta: 0.01(33)$
}
 \\
\midrule
0.2 & \shortstack{
$a: -29.74(17)$\\
$b: -2.17(17)$\\
$c: 13.24(17)$\\
$d: -0.11(17)$\\
$\delta: 0.01(17)$
} & \shortstack{
$a: -29.05(20)$\\
$b: -2.19(18)$\\
$c: 12.81(18)$\\
$d: -0.10(18)$\\
$\delta: 0.01(18)$
} & \shortstack{
$a: -24.42(40)$\\
$b: -0.73(36)$\\
$c: 9.63(40)$\\
$d: -0.12(36)$\\
$\delta: 0.00(36)$
} \\
\bottomrule
\end{tabular}
}
\end{tabular}
\end{center}
\end{table}

\begin{table}
\begin{center}
\caption{Estimated parameters of the linear model \ref{eq:linearPDE_app} to which spatio-temporal data of fluctuations were fitted. This table denotes the values for the same settings as shown in Figure \ref{Fig:Noise_level_correlation_Tipping}, i.e. Experiment 1 from Table \ref{tab:experiments_final} with $\delta=0.5$ (Tipping case).} 
\label{Table:Noise_level_correlation_Tipping}
\begin{tabular}{@{}c@{\hskip 0.5em}c@{}}
\raisebox{-3.5em}{\rotatebox{90}{\textbf{Correlation Length}}} & \shortstack{\textbf{Noise Level} \\
\begin{tabular}{|c|c|c|c|}
\toprule
\textbf{} & 0.1 & 1 & 10 \\
\midrule
\textbf{\rotatebox{90}{White Noise}} &
\shortstack{
$a: -29.79(36)$\\
$b: -2.21(36)$\\
$c: 13.27(36)$\\
$d: -0.12(36)$\\
$\delta: 0.50(36)$
} & \shortstack{
$a: -29.73(36)$\\
$b: -2.21(36)$\\
$c: 13.23(36)$\\
$d: -0.12(36)$\\
$\delta: 0.50(36)$
} & \shortstack{
$a: -25.93(40)$\\
$b: -2.20(38)$\\
$c: 10.15(38)$\\
$d: -0.09(38)$\\
$\delta: 0.50(38)$
} \\
\midrule
0.1 & \shortstack{
$a: -29.78(37)$\\
$b: -2.20(37)$\\
$c: 13.25(37)$\\
$d: -0.14(37)$\\
$\delta: 0.50(37)$
} & \shortstack{
$a: -29.64(37)$\\
$b: -2.21(37)$\\
$c: 13.16(37)$\\
$d: -0.13(37)$\\
$\delta: 0.50(37)$
} & \shortstack{
$a: -22.10(40)$\\
$b: -1.82(39)$\\
$c: 7.84(39)$\\
$d: -0.25(39)$\\
$\delta: 0.50(39)$
} \\
\midrule
0.2 & \shortstack{
$a: -29.74(38)$\\
$b: -2.20(38)$\\
$c: 13.24(38)$\\
$d: -0.18(38)$\\
$\delta: 0.49(38)$
} & \shortstack{
$a: -29.45(40)$\\
$b: -2.21(38)$\\
$c: 13.08(38)$\\
$d: -0.17(38)$\\
$\delta: 0.49(38)$
} & \shortstack{
$a: -19.51(40)$\\
$b: -1.41(44)$\\
$c: 6.79(40)$\\
$d: -0.35(44)$\\
$\delta: 0.48(44)$
} \\
\bottomrule
\end{tabular}
}
\end{tabular}
\end{center}
\end{table}

\begin{figure}[h!]
  \centering
  \includegraphics[width=\textwidth]{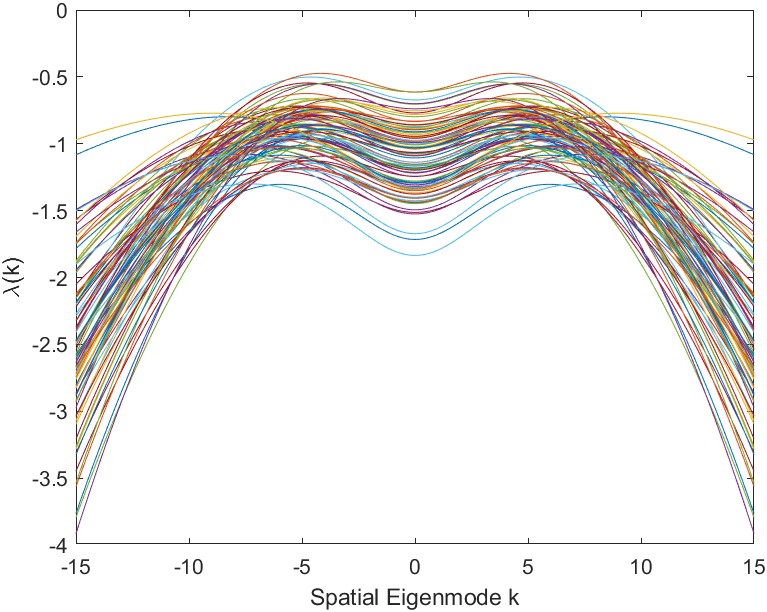}
  \caption{All individual dispersion relations found for Experiment 1 with Noise Strength $0.1$, correlation length $l_c = 0.2$ and $\delta=0.01$ (Turing case). Figure \ref{Fig:Noise_level_correlation_Turing} reports on the statistics of these. Noteworthy here is that some of these estimated dispersion relation have different forms with peaks for different spatial eigenmodes $k_*$, which might explain the spread of the eigenmodes found in the stastics.}
  \label{fig:correlation_0.2}
\end{figure}

\begin{figure}[h!]
  \centering
  \includegraphics[width=0.9\textwidth]{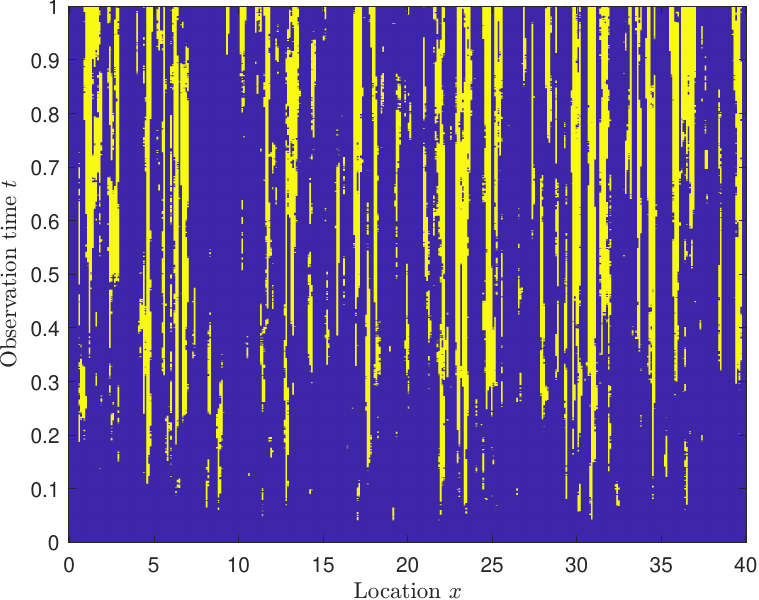}
  \caption{Heatmap of a random run from the Klausmeier model (Eq.~\ref{eq:Extended_klausmeier_2}) with noise strength 10 and baseline values from Table~\ref{tab:experiments_final}. Blue indicates regions where the vegetation variable $v$ remains above the unstable state threshold ($v>0.0846469$), while yellow marks regions where it falls below this threshold ($v<0.0846469$). This demonstrates that strong noise can locally push the vegetation towards different attractors.}
  \label{fig:Heatmap_effect_Noise}
\end{figure}

\clearpage
\subsection{Additional Material for Numerical Experiments 2 \& 3}\label{Appendix:Experiment_23}

\begin{table}[H]
\begin{center}
\caption{Estimated parameters of the linear model \ref{eq:linearPDE_app} to which spatio-temporal data of fluctuations were fitted. This table denotes the values for the same settings as shown in Figure \ref{Fig:sampling_Turing}, i.e. Experiment 2 from Table \ref{tab:experiments_final} with $\delta=0.01$ (Turing case).} 
\label{Table:Sampling_Turing}
\begin{tabular}{@{}c@{\hskip 0.5em}c@{}}
\raisebox{-3.5em}{\rotatebox{90}{\textbf{Observation Time}}} & \shortstack{\textbf{Temporal sampling} \\
\begin{tabular}{|c|c|c|c|}
\toprule
\textbf{} & $100\%$ & $50\%$ & $25\%$ \\
\midrule
0.2 & \shortstack{
$a: -29.94(70)$\\
$b: -2.16(66)$\\
$c: 13.28(70)$\\
$d: -0.10(66)$\\
$\delta: 0.01(66)$
} & \shortstack{
$a: -29.52(70)$\\
$b: -2.15(66)$\\
$c: 13.24(70)$\\
$d: -0.10(66)$\\
$\delta: 0.01(66)$
} & \shortstack{
$a: -28.74(70)$\\
$b: -2.13(66)$\\
$c: 13.16(70)$\\
$d: -0.11(66)$\\
$\delta: 0.01(66)$
} \\
\midrule
1 & \shortstack{
$a: -29.38(17)$\\
$b: -2.17(17)$\\
$c: 12.93(17)$\\
$d: -0.08(17)$\\
$\delta: 0.01(17)$
} & \shortstack{
$a: -29.02(17)$\\
$b: -2.16(17)$\\
$c: 12.89(17)$\\
$d: -0.09(17)$\\
$\delta: 0.01(17)$
} & \shortstack{
$a: -28.30(17)$\\
$b: -2.13(17)$\\
$c: 12.81(17)$\\
$d: -0.09(17)$\\
$\delta: 0.01(17)$
} \\
\midrule
5 & \shortstack{
$a: -28.931(60)$\\
$b: -2.196(59)$\\
$c: 12.793(60)$\\
$d: -0.043(59)$\\
$\delta: 0.010(59)$
} & \shortstack{
$a: -28.590(60)$\\
$b: -2.183(59)$\\
$c: 12.751(60)$\\
$d: -0.045(59)$\\
$\delta: 0.010(59)$
} & \shortstack{
$a: -27.916(60)$\\
$b: -2.157(59)$\\
$c: 12.669(60)$\\
$d: -0.049(59)$\\
$\delta: 0.010(59)$
} \\
\bottomrule
\end{tabular}
}
\end{tabular}
\end{center}
\end{table}

\begin{table}
\begin{center}
\caption{Estimated parameters of the linear model \ref{eq:linearPDE_app} to which spatio-temporal data of fluctuations were fitted. This table denotes the values for the same settings as shown in Figure \ref{Fig:sampling_Tipping}, i.e. Experiment 2 from Table \ref{tab:experiments_final} with $\delta=0.5$ (Tipping case).} 
\begin{tabular}{@{}c@{\hskip 0.5em}c@{}}
\raisebox{-3.5em}{\rotatebox{90}{\textbf{Observation Time}}} & \shortstack{\textbf{Temporal sampling} \\
\begin{tabular}{|c|c|c|c|}
\toprule
\textbf{} & $100\%$ & $50\%$ & $25\%$ \\
\midrule
0.2 & \shortstack{
$a: -29.94(70)$\\
$b: -2.16(66)$\\
$c: 13.28(70)$\\
$d: -0.10(66)$\\
$\delta: 0.01(66)$
} & \shortstack{
$a: -29.52(70)$\\
$b: -2.15(66)$\\
$c: 13.24(70)$\\
$d: -0.10(66)$\\
$\delta: 0.01(66)$
} & \shortstack{
$a: -28.74(70)$\\
$b: -2.13(66)$\\
$c: 13.16(70)$\\
$d: -0.11(66)$\\
$\delta: 0.01(66)$
} \\
\midrule
1 & \shortstack{
$a: -29.38(17)$\\
$b: -2.17(17)$\\
$c: 12.93(17)$\\
$d: -0.08(17)$\\
$\delta: 0.01(17)$
} & \shortstack{
$a: -29.02(17)$\\
$b: -2.16(17)$\\
$c: 12.89(17)$\\
$d: -0.09(17)$\\
$\delta: 0.01(17)$
} & \shortstack{
$a: -28.30(17)$\\
$b: -2.13(17)$\\
$c: 12.81(17)$\\
$d: -0.09(17)$\\
$\delta: 0.01(17)$
} \\
\midrule
5 & \shortstack{
$a: -28.931(60)$\\
$b: -2.196(59)$\\
$c: 12.793(60)$\\
$d: -0.043(59)$\\
$\delta: 0.010(59)$
} & \shortstack{
$a: -28.590(60)$\\
$b: -2.183(59)$\\
$c: 12.751(60)$\\
$d: -0.045(59)$\\
$\delta: 0.010(59)$
} & \shortstack{
$a: -27.916(60)$\\
$b: -2.157(59)$\\
$c: 12.669(60)$\\
$d: -0.049(59)$\\
$\delta: 0.010(59)$
} \\
\bottomrule
\end{tabular}
}
\end{tabular}
\end{center}
\label{Table:sampling_Tipping}
\end{table}

\begin{table}
\begin{center}
\caption{Estimated parameters of the linear model \ref{eq:linearPDE_app} to which spatio-temporal data of fluctuations were fitted. This table denotes the values for the same settings as shown in Figure \ref{Fig:spatial_sampling_corr_02}, i.e. Experiment 3 from Table \ref{tab:experiments_final}.} 
\begin{tabular}{@{}c@{\hskip 0.5em}c@{}}
\raisebox{-5em}{\rotatebox{90}{\textbf{Diffusion parameter $\delta$}}} & \shortstack{\textbf{Spatial sampling} \\
\begin{tabular}{|c|c|c|c|}
\toprule
\textbf{} & $100\%$ & $50\%$ & $25\%$ \\
\midrule
0.01 & \shortstack{
$a: -29.65(15)$\\
$b: -2.17(15)$\\
$c: 13.02(15)$\\
$d: -0.07(15)$\\
$\delta: 0.01(15)$
} & \shortstack{
$a: -30.21(15)$\\
$b: -2.18(15)$\\
$c: 13.03(15)$\\
$d: -0.05(15)$\\
$\delta: 0.01(15)$
} & \shortstack{
$a: -32.23(20)$\\
$b: -2.20(16)$\\
$c: 13.04(16)$\\
$d: -0.00(16)$\\
$\delta: 0.01(16)$
} \\
\midrule
0.5 & \shortstack{
$a: -29.91(31)$\\
$b: -2.15(31)$\\
$c: 13.26(31)$\\
$d: -0.10(31)$\\
$\delta: 0.50(31)$
} & \shortstack{
$a: -30.47(30)$\\
$b: -2.17(30)$\\
$c: 13.29(30)$\\
$d: 0.02(30)$\\
$\delta: 0.52(30)$
} & \shortstack{
$a: -32.54(26)$\\
$b: -2.22(26)$\\
$c: 13.41(26)$\\
$d: 0.47(26)$\\
$\delta: 0.58(26)$
} \\
\bottomrule
\end{tabular}
}
\end{tabular}
\label{Table:spatial_sampling_corr_02}
\end{center}
\end{table}

\begin{figure}[h!]
\includegraphics[width=\textwidth]{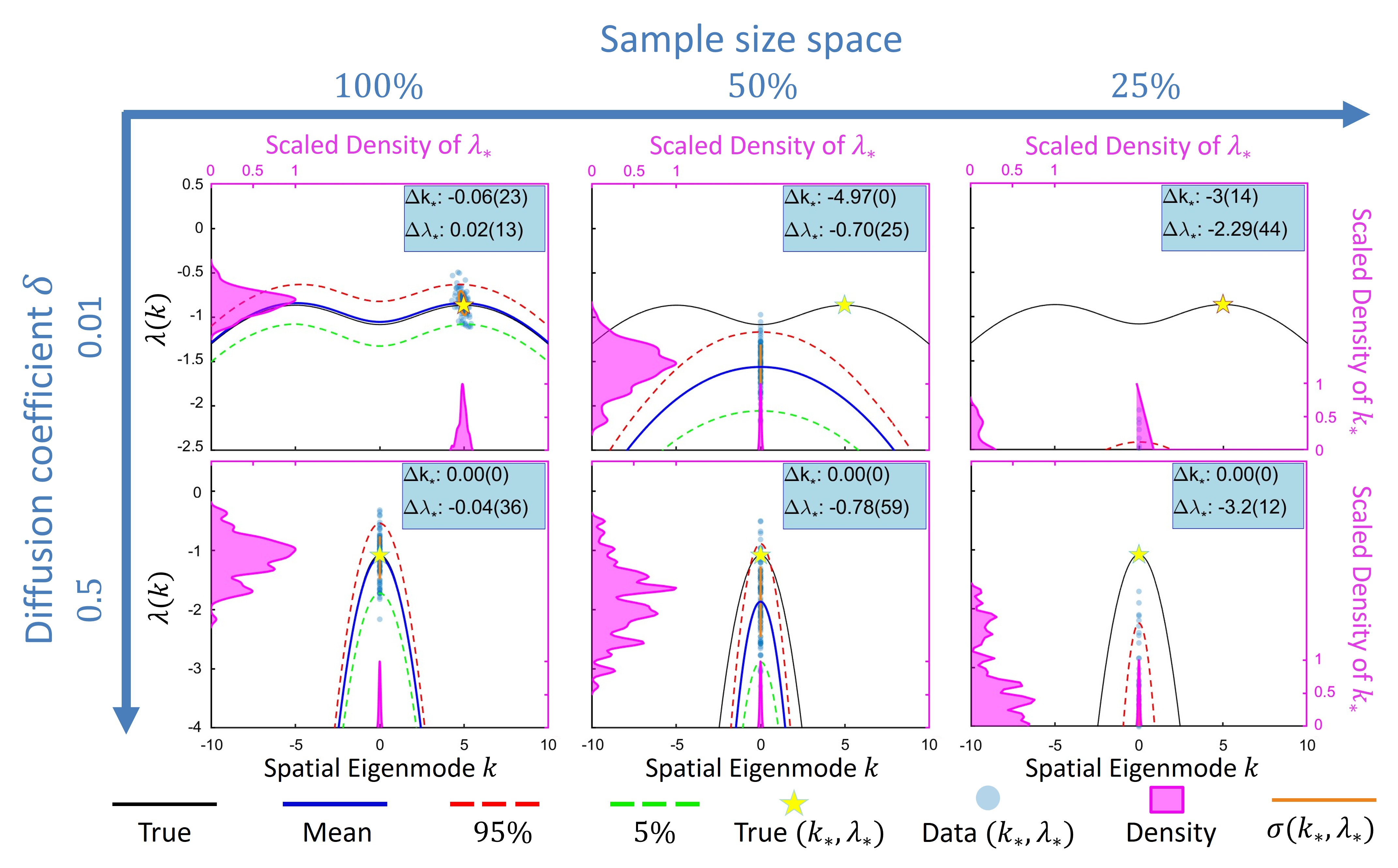}
\caption{Estimated dispersion relations for different spatial sampling rates for $\delta = 0.01$ (Turing case; top row) and $\delta = 0.5$ (saddle-node case; bottom row). Data is generated using Experiment 3 settings in Table~\ref{tab:experiments_final}, except the added noise is now spatially uncorrelated white noise $(l_c=0)$. See the caption of Figure \ref{Fig:Noise_level_correlation_Turing} for the details of the depicted lines, areas, circles and insets.}
\label{Fig:spatial_sampling_white}
\end{figure}

\begin{table}
\begin{center}
\caption{Estimated parameters of the linear model \ref{eq:linearPDE_app} to which spatio-temporal data of fluctuations were fitted. This table denotes the values for the same settings as shown in Figure \ref{Fig:spatial_sampling_white}, i.e. Experiment 3 from Table \ref{tab:experiments_final}, except now $l_c=0$.} 
\label{Table:spatial_sampling_white}
\begin{tabular}{@{}c@{\hskip 0.5em}c@{}}
\raisebox{-5em}{\rotatebox{90}{\textbf{Diffusion parameter $\delta$}}} & \shortstack{\textbf{Spatial sampling} \\
\begin{tabular}{|c|c|c|c|}
\toprule
\textbf{} & $100\%$ & $50\%$ & $25\%$ \\
\midrule
0.01 & \shortstack{
$a: -29.72(16)$\\
$b: -2.16(16)$\\
$c: 13.08(16)$\\
$d: -0.10(16)$\\
$\delta: 0.01(16)$
} & \shortstack{
$a: -79.64(30)$\\
$b: -2.57(27)$\\
$c: 12.52(27)$\\
$d: -1.03(27)$\\
$\delta: 0.01(27)$
} & \shortstack{
$a: -101.19(30)$\\
$b: -2.61(28)$\\
$c: 12.39(30)$\\
$d: -1.40(28)$\\
$\delta: 0.01(28)$
} \\
\midrule
0.5 & \shortstack{
$a: -29.85(36)$\\
$b: -2.18(36)$\\
$c: 13.25(36)$\\
$d: -0.11(36)$\\
$\delta: 0.50(36)$
} & \shortstack{
$a: -80.79(60)$\\
$b: -3.41(60)$\\
$c: 13.00(60)$\\
$d: -1.34(60)$\\
$\delta: 0.98(60)$
} & \shortstack{
$a: -102.7(13)$\\
$b: -4.0(13)$\\
$c: 12.0(13)$\\
$d: -4.2(13)$\\
$\delta: 1.8(13)$
} \\
\bottomrule
\end{tabular}
}
\end{tabular}
\end{center}
\end{table}
\clearpage
\subsection{Additional Material for Numerical Experiment 4}\label{Appendix:Experiment_4}

\begin{figure}[H]
    \centering
    \includegraphics[width=0.7\linewidth]{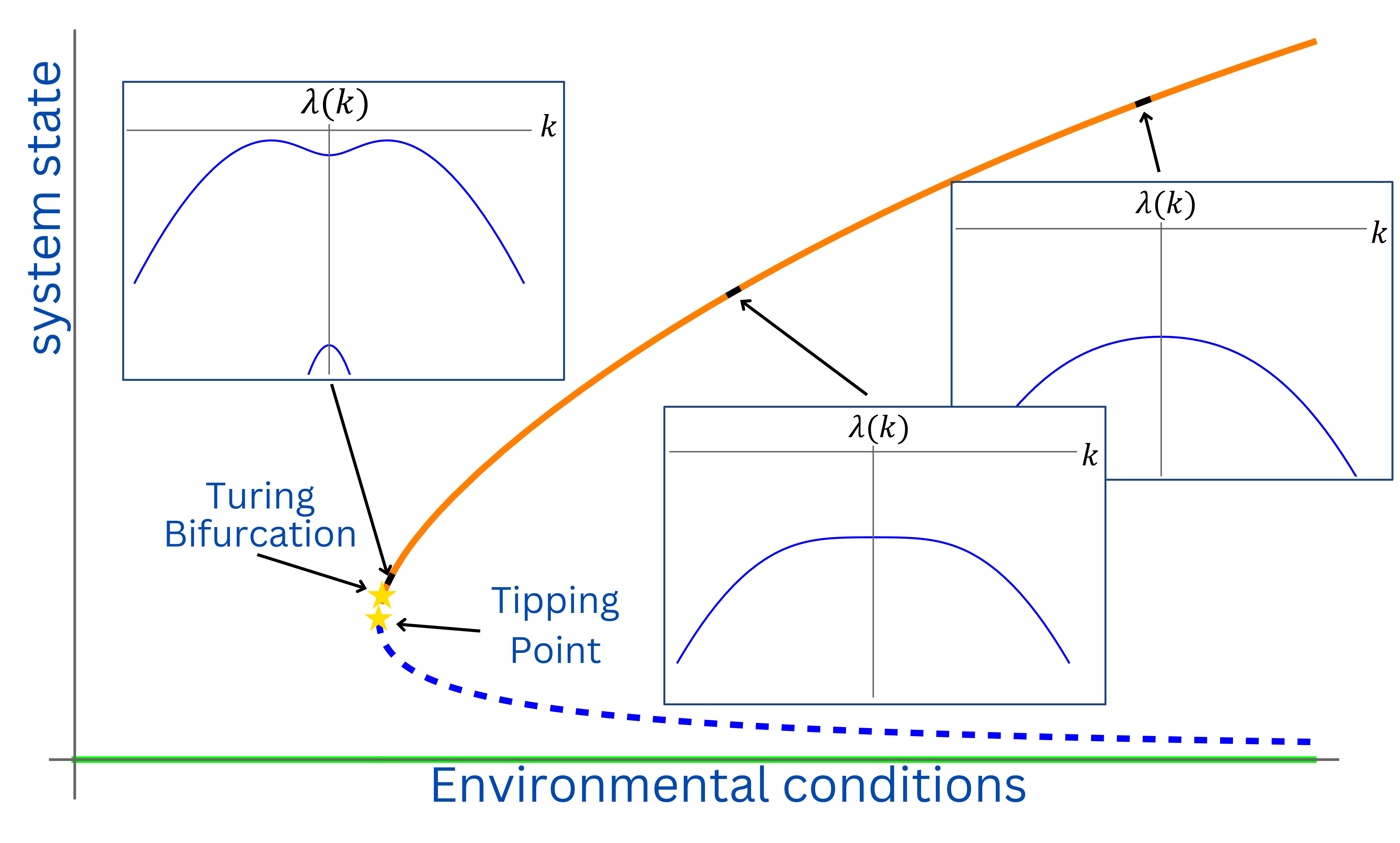}
    \caption{Bifurcation diagram including insets with dispersion relations for the case $\delta = 0.1$. For these parameter settings, the orange line does destabilise via a Turing bifurcation, but this bifurcation is located very close to the saddle-node bifurcation ($p_\text{T} \approx 1.11445$ and $p_\text{SN} \approx 1.10499$). Hence the dispersion relation only shows a peak for $k_* \neq 0$ very close to the Turing bifurcation, and it can appear very flat around $k = 0$. Note that Figure~\ref{Fig:Early_Warning_Overview} shows similar figures for cases $\delta = 0.01$ (Turing case) and $\delta = 0.5$ (Tipping case).} 
    \label{fig:Early_warning_shift}
\end{figure}

\begin{table}[H]
\begin{center}
\caption{Estimated parameters of the linear model \ref{eq:linearPDE_app} to which spatio-temporal data of fluctuations were fitted. This table denotes the values for the same settings as shown in Figure \ref{Fig:parameter_tests}, i.e. Experiment 4 from Table \ref{tab:experiments_final}.} 
\begin{tabular}{@{}c@{\hskip 0.5em}c@{}}
\raisebox{-6.5em}{\rotatebox{90}{\textbf{Environmental condition $p$}}} & \shortstack{\textbf{Diffusion coefficient $\delta$} \\
\begin{tabular}{|c|c|c|c|}
\toprule
\textbf{} & 0.1 & 1 &  10 \\
\midrule
\textbf{\rotatebox{90}{20}} &
 \shortstack{
$a: -64.16(24)$\\
$b: -4.95(24)$\\
$c: 12.62(20)$\\
$d: -1.47(24)$\\
$\delta: 0.01(24)$
} & \shortstack{
$a: -64.32(35)$\\
$b: -4.96(35)$\\
$c: 12.68(30)$\\
$d: -1.48(35)$\\
$\delta: 0.10(35)$
} & \shortstack{
$a: -64.41(44)$\\
$b: -4.98(44)$\\
$c: 12.71(40)$\\
$d: -1.51(44)$\\
$\delta: 0.50(44)$
} \\
\midrule
6 & \shortstack{
$a: -29.38(17)$\\
$b: -2.17(17)$\\
$c: 12.93(17)$\\
$d: -0.08(17)$\\
$\delta: 0.01(17)$
} & \shortstack{
$a: -29.54(28)$\\
$b: -2.18(28)$\\
$c: 13.09(28)$\\
$d: -0.10(28)$\\
$\delta: 0.10(28)$
} & \shortstack{
$a: -29.64(37)$\\
$b: -2.21(37)$\\
$c: 13.16(37)$\\
$d: -0.13(37)$\\
$\delta: 0.50(37)$
} \\
\midrule
2 & \shortstack{
$a: -7.60(14)$\\
$b: -1.30(14)$\\
$c: 4.68(14)$\\
$d: 0.26(14)$\\
$\delta: 0.01(14)$
} & \shortstack{
$a: -7.54(28)$\\
$b: -1.35(28)$\\
$c: 4.73(28)$\\
$d: 0.28(28)$\\
$\delta: 0.10(28)$
} & \shortstack{
$a: -7.57(43)$\\
$b: -1.39(43)$\\
$c: 4.78(43)$\\
$d: 0.26(43)$\\
$\delta: 0.50(43)$
} \\
\bottomrule
\end{tabular}
}
\end{tabular}
\end{center}
\label{Table:parameters_testss}
\end{table}

\clearpage
\subsection{Additional Material for Numerical Experiment 5}\label{Appendix:Experiment_5}

\begin{table}[H]
\begin{center}
\caption{Estimated parameters of the linear model \ref{eq:linearPDE_app} to which spatio-temporal data of fluctuations were fitted. This table denotes the values for the same settings as shown in Figure \ref{Fig:parameter_time_variatiations}, i.e. Experiment 5 from Table \ref{tab:experiments_final}.} 
\begin{tabular}{@{}c@{\hskip 0.5em}c@{}}
\raisebox{-8.5em}{\rotatebox{90}{\textbf{Environmental condition interval $p$}}} & \shortstack{\textbf{Diffusion coefficient $\delta$} \\
\begin{tabular}{|c|c|c|c|}
\toprule
\textbf{} & 0.1 & 1 & 10 \\
\midrule
\textbf{\rotatebox{90}{[20,18]}} &
 \shortstack{
$a: -63.53(24)$\\
$b: -4.71(24)$\\
$c: 12.73(20)$\\
$d: -1.39(24)$\\
$\delta: 0.01(24)$
} & \shortstack{
$a: -63.69(35)$\\
$b: -4.72(35)$\\
$c: 12.80(30)$\\
$d: -1.40(35)$\\
$\delta: 0.10(35)$
} & \shortstack{
$a: -63.78(44)$\\
$b: -4.75(44)$\\
$c: 12.83(40)$\\
$d: -1.43(44)$\\
$\delta: 0.50(44)$
} \\
\midrule
\textbf{\rotatebox{90}{[8,6]}} & \shortstack{
$a: -29.38(17)$\\
$b: -2.17(17)$\\
$c: 12.93(17)$\\
$d: -0.08(17)$\\
$\delta: 0.01(17)$
} & \shortstack{
$a: -36.67(29)$\\
$b: -2.27(29)$\\
$c: 14.14(29)$\\
$d: -0.29(29)$\\
$\delta: 0.10(29)$
} & \shortstack{
$a: -36.76(37)$\\
$b: -2.30(37)$\\
$c: 14.20(37)$\\
$d: -0.32(37)$\\
$\delta: 0.50(37)$
} \\
\midrule
\textbf{\rotatebox{90}{[4,2]}} & \shortstack{
$a: -18.12(16)$\\
$b: -1.29(16)$\\
$c: 9.73(16)$\\
$d: -0.03(16)$\\
$\delta: 0.01(16)$
} & \shortstack{
$a: -18.22(28)$\\
$b: -1.34(28)$\\
$c: 9.87(28)$\\
$d: -0.04(28)$\\
$\delta: 0.10(28)$
} & \shortstack{
$a: -18.29(36)$\\
$b: -1.37(36)$\\
$c: 9.95(36)$\\
$d: -0.06(36)$\\
$\delta: 0.50(36)$
}\\
\bottomrule
\end{tabular}
}
\end{tabular}
\end{center}
\label{Table:parameter_time_variations}
\end{table}

\end{document}